\documentclass[a4paper,11pt]{article}

\usepackage{xcolor}
\usepackage{graphicx}

\usepackage{amscd}
\usepackage{amsmath}
\usepackage{ascmac}
\usepackage{cases}
\usepackage{amssymb}
\usepackage{amsfonts}
\usepackage{amsthm}
\usepackage{bm}

\usepackage[all,pdf]{xy}

\usepackage{enumerate}
\usepackage{wrapfig}

\usepackage{tikz}
\usepackage{tikz-cd}

\definecolor{Mylinkclr}{HTML}{666666} 

\usepackage[
    pagebackref,
    colorlinks=true,
    linkcolor=Mylinkclr,
    citecolor=Mylinkclr,
    linktoc=all
]{hyperref}

\usepackage{zref-clever}
\zcsetup{cap}

\zcRefTypeSetup{equation}{
    Name-sg=,
    name-sg=,
    Name-pl=,
    name-pl=,
    Name-sg-ab=,
    name-sg-ab=,
    Name-pl-ab=,
    name-pl-ab=
}

\theoremstyle{plain}

\newtheorem{Thm}{{\textbf{Theorem}}}[section]
\AddToHook{env/Thm/begin}{%
    \zcsetup{countertype={Thm=Thm}}%
}
\zcRefTypeSetup{Thm}{Name-sg=}

\newtheorem{Lem}[Thm]{{\textbf{Lemma}}}
\AddToHook{env/Lem/begin}{%
    \zcsetup{countertype={Thm=Lem}}%
}
\zcRefTypeSetup{Lem}{Name-sg=}

\newtheorem{Prop}[Thm]{{\textbf{Proposition}}}
\AddToHook{env/Prop/begin}{%
    \zcsetup{countertype={Thm=Prop}}%
}
\zcRefTypeSetup{Prop}{Name-sg=}

\newtheorem{Cor}[Thm]{{\textbf{Corollary}}}
\AddToHook{env/Cor/begin}{%
    \zcsetup{countertype={Thm=Cor}}%
}
\zcRefTypeSetup{Cor}{Name-sg=}

\newtheorem{Conj}[Thm]{{\textbf{Conjecture}}}
\AddToHook{env/Conj/begin}{%
    \zcsetup{countertype={Thm=Conj}}%
}
\zcRefTypeSetup{Conj}{Name-sg=}

\newtheorem{Def}[Thm]{{\textbf{Definition}}}
\AddToHook{env/Def/begin}{%
    \zcsetup{countertype={Thm=Def}}%
}
\zcRefTypeSetup{Def}{Name-sg=}

\newtheorem{Rem}[Thm]{{\textbf{Remark}}}
\AddToHook{env/Rem/begin}{%
    \zcsetup{countertype={Thm=Rem}}%
}
\zcRefTypeSetup{Rem}{Name-sg=}

\AddToHook{env/Ex/begin}{%
    \zcsetup{countertype={Thm=Ex}}%
}
\zcRefTypeSetup{Ex}{Name-sg=}

\newcommand{\e}{\epsilon}

\newcommand{\ba}{{\bf a}}

\newcommand{\bnu}{{\boldsymbol \nu}}

\newcommand{\bl}{{\bm l}}

\newcommand{\bC}{{\mathbb C}}       
\newcommand{\bP}{{\mathbb P}}      
\newcommand{\bZ}{{\mathbb Z}}

\newcommand{\cA}{{\mathcal A}}
\newcommand{\cG}{{\mathcal G}}
\newcommand{\cO}{{\mathcal O}}

\newcommand{\cM}{{\mathcal M}}

\newcommand{\cD}{{\mathcal D}}

\newcommand{\cJ}{{\mathcal J}}

\newcommand{\cF}{{\mathcal F}}
\newcommand{\cE}{{\mathcal E}}

\newcommand{\cQ}{{\mathcal Q}}

\newcommand{\cL}{{\mathcal L}}
\newcommand{\cY}{{\mathcal Y}}

\newcommand{\Aut}{\mathop{\rm Aut}\nolimits}

\newcommand{\tr}{\mathop{\rm tr}\nolimits}

\newcommand\End{\mathop{\rm End}\nolimits}

\newcommand\Pic{\mathop{\rm Pic}\nolimits}
\newcommand\Spec{\mathop{\rm Spec}\nolimits}

\newcommand\Ker{\mathop{\rm Ker}\nolimits}

\newcommand{\res}{\mathop{\sf res}\nolimits}

\newcommand{\Elm}{\mathop{\sf Elm}\nolimits}

\newcommand{\App}{ \mathop{\sf App}\nolimits}
\newcommand{\norm}{ \mathop{\rm norm}\nolimits}
\newcommand{\Sym}{ \mathop{\rm Sym}\nolimits}

\newcommand{\Supp}{\mathop{\rm Supp}\nolimits}
\newcommand{\id}{\mathop{\rm id}\nolimits}

\newcommand{\Coh}{\mathop{\rm Coh}\nolimits}
\newcommand{\Stab}{\mathop{\rm Stab}\nolimits}
\newcommand{\Div}{\mathop{\rm div}\nolimits}
\newcommand{\Odd}{\mathop{\rm Odd}\nolimits}
\newcommand{\Oka}{\mathop{\rm Oka}\nolimits}
\newcommand{\Par}{\mathop{\rm Par}\nolimits}

\newcommand\ra{\rightarrow}

\newcommand{\Bun}{ \mathop{\sf Bun}\nolimits}

\newcommand{\Address}{{
  \bigskip
  \footnotesize
  
  Yuki Matsubara, \textsc{Centre for Quantum Mathematics, University of Southern Denmark, Campusvej 55, DK-5230, Odense. M, Denmark}\par\nopagebreak
  \textit{E-mail address}, Yuki Matsubara: \texttt{matsubara@imada.sdu.dk}
}}

\begin{document}
\title{Cohomology of vector bundles on the moduli space of parabolic connections on $\bP^1$ minus $5$ points}
\author{Yuki Matsubara}
\date{}
\maketitle

\begin{abstract}
We study the moduli space of parabolic connections of rank two on the complex projective line $\mathbb{P}^1$ minus five points with fixed spectral data. 
We prove a cohomological vanishing theorem for tensor products of four natural vector bundles on the moduli stack of $\bnu$-${\mathfrak sl}_2$-parabolic connections.
We also compute the cohomology of the structure sheaf and give an explicit geometric description of a compactification of this moduli space, using $\epsilon$-connections and the birational geometry of $\bP^2 \times (\bP^2)^{\vee}$.
As an application, we construct a universal family of twisted $D$-modules whose fibers satisfy an orthogonality property over a dense open subset.
These results are motivated by Arinkin's approach to the tamely ramified geometric Langlands correspondence.

\end{abstract}

\tableofcontents

\section{Introduction}
\subsection{Overview}
In this paper, we study moduli spaces of rank-two parabolic connections on the complex projective line minus $n$ points with fixed spectral data.
More precisely, we consider logarithmic connections on $\bP^1$ with poles at $n$ marked points $\{ t_1, \dots, t_n \}$, fixed determinant, and generic local exponents.
These moduli spaces appear naturally as spaces of initial conditions for Garnier systems.
The case $n = 4$ is closely related to the Okamoto initial condition space of Painlev\'{e} V\hspace{-1.2pt}I equation \cite{Oka79}, and has been studied extensively by Arinkin--Lysenko \cite{AL97}, Inaba--Iwasaki--Saito \cite{IISii}, and others.

The main purpose of this paper is to study the geometry and cohomology of these moduli spaces in the case $n = 5$.
In this case, the corresponding moduli space has complex dimension four.
We compactify the moduli space by allowing $\epsilon$-connections.
The boundary locus at $\epsilon = 0$ is described in terms of parabolic Higgs bundles, spectral curves, and twisted compactified Jacobians.
In the five points case, we give an explicit birational description of the compactification using the apparent map and the forgetful map to the moduli space of quasi-parabolic bundles.
This relates the compactification to the birational geometry of $\bP^2 \times (\bP^2)^{\vee}$, extending the result given by Loray--Saito \cite{LS}.

Our main result regarding cohomology is a vanishing theorem for natural vector bundles on the moduli stack $\cM$ of $\bnu$-${\mathfrak sl}_2$-parabolic connections.
For a point $x \in \bP^1$, let $\xi_x$ be the rank two vector bundle on $\cM$ obtained by evaluating the universal bundle at $x$.
When $n = 5$, we prove that, for four pairwise distinct points $x_1, \dots, x_4 \in \bP^1$, one has
\begin{equation*}
 H^i(\cM, \xi_{x_1} \otimes \xi_{x_2} \otimes \xi_{x_3} \otimes \xi_{x_4}) = 0
\end{equation*}
for all $i \geq 0$.
We also compute the cohomology of the structure sheaf of $\cM$.
These results give the cohomological orthogonality statements needed in the construction considered below.

Finally, we explain the relation with Arinkin's approach to the tamely ramified geometric Langlands correspondence.
Let $P^{\vee}$ be the non-separated scheme obtained by gluing two copies of $\bP^{n-3}$ along the complement of the hyperplanes $Z_i$ determined by the marked points.
For every $n \geq 5$, we construct a universal family $\xi_{\bnu}$ of twisted $D$-modules on $\cM \times P^{\vee}$.
In the case $n = 5$, the vanishing theorem above implies that the fibers of this family satisfy an orthogonality property over a dense open subset of $P^{\vee}$.
This provides evidence for the expected Fourier--Mukai transform with the kernel $\xi_{\bnu}$, motivated by Arinkin's work in the $n = 4$ case \cite{A01}.

\subsection{Summary of main results}
In order to state the main theorem precisely, we will introduce some notations.
Fix $t_1, \dots, t_n \in \bP^1$ and $\nu_1, \dots, \nu_n \in \bC$ such that $t_i \neq t_j$ for $i \neq j$, $n \geq 4$, $2\nu_i \not\in \mathbb{Z}$, and
\begin{equation}
 \sum_{i = 1}^n \epsilon_i \nu_i \not\in \bZ
\end{equation}
for any $\epsilon_i \in \mu_2 := \{1, -1\}$.

\begin{Def}
A {\rm\bf{$\bnu$-$\mathfrak{s}l_2$-parabolic connection}} is a triple $(L, \nabla, \varphi)$ such that
\begin{itemize}
 \item[(1)] $L$ is a rank $2$ vector bundle on $\bP^1$,
 \item[(2)] $\nabla: L \ra L \otimes \Omega^1_{\bP^1}(D)$ is a connection, where $D := t_1 + \cdots + t_n$,
 \item[(3)] $\varphi : \bigwedge^2 L \xrightarrow{\sim} \cO_{\bP^1}$ is an isomorphism which satisfies
 \begin{equation*}
  \varphi (\nabla s_1 \wedge s_2) + \varphi(s_1 \wedge \nabla s_2) = d(\varphi(s_1 \wedge s_2))
 \end{equation*}
 for $s_1, s_2 \in L$,
 \item[(4)] the residue $\res_{t_i}(\nabla)$ of the connection $\nabla$ at $t_i$ has eigenvalues $\{\nu_i, -\nu_i \}$ for each $i$ $(1 \leq i \leq n)$.
\end{itemize}
\end{Def}

Let $\cM$ be the moduli stack of $\bnu$-$\mathfrak{s}l_2$-parabolic connections, and $M$ be the corresponding coarse moduli space.
It is known that $M$ is a smooth irreducible separated quasi-projective scheme of dimension $2(n-3)$, and $\cM$ is a $\mu_2$-gerbe over $M$ (\cite{AL97, IISi}).

Let $Z_i^{\pm}$ be the pre-images of $Z_i \subset \bP^{n-3}$ along $p \colon P^{\vee} \ra \bP^{n-3}$, and
$\bnu := \sum_{i = 1}^n \nu_i ([Z_i^+] - [Z_i^-]) \in \mathrm{div}(P^{\vee}) \otimes_{\bZ} \bC$, where $\mathrm{div}(P^{\vee})$ is the group of divisors on $P^{\vee}$.
Let $D_{\bnu}$ denote the twisted differential operator (TDO) ring corresponding to $\bnu$ over $P^{\vee}$.

Firstly, we will explain our construction of a universal $D_{\bnu}$-module $\xi_{\bnu}$ over $\cM \times P^{\vee}$.
For any connection $\mathbb{L} = (L, \nabla, \varphi) \in \mathcal{M}$, its symmetric product $\Sym^{n-3}(\mathbb{L})$ gives a connection on $\mathbb{P}^{n-3}$.
More precisely, it is the symmetric part of the push-forward of $\mathbb{L}^{\boxtimes n-3}$ along the map $\Sym \colon (\bP^1)^{n-3} \ra \bP^{n-3}$, that is, $\Sym^{n-3}(\mathbb{L}) := (\Sym_*(\mathbb{L}^{\boxtimes (n-3)}))^{\mathfrak{S}_{n-3}}$. 
This connection has singularities along the divisors $Z_i$ $(i = 1, \dots, n)$, as well as along the discriminant divisor $\Delta \subset \bP^{n-3}$.
The divisors $Z_i$ cross normally, and the singularity along $Z_i$ has residue with eigenvalues $\{\pm \nu_i \}$, each with multiplicity $2^{n-4}$.
Let us construct the $D_{\bnu}$-module $ j_{!*}(\Sym^{n-3}(\mathbb{L})|_U)$ with $j \colon U := \bP^{n-3} \setminus \cup_{i = 1}^n Z_i \hookrightarrow P^{\vee}$.
This construction still makes sense for a family of connections.
Let us apply it to the universal family of connections, and get a $\cM$-family $\xi_{\bnu}$ of $D_{\bnu}$-modules over $\cM \times P^{\vee}$.

Suppose $n=5$. 
The main theorem in this paper is as follows:
For $x \in \bP^1$ let $\xi_x$ be the bundle on $\cM$ whose fiber at $(L, \nabla, \varphi)$ is $L_x$.

\begin{Thm}\label{main}
Suppose $x_1,\dots, x_4 \in \bP^1$ and $x_i \neq x_j$ for $i \neq j$. Then
\begin{equation*}
 H^i(\cM, \xi_{x_1}\otimes \xi_{x_2} \otimes \xi_{x_3} \otimes \xi_{x_4}) = 0
\end{equation*}
 for any $i \geq 0$.
\end{Thm}
It is predicted that for general $n \geq 4$, the similar statement is also true with $2(n-3)$ points on $\bP^1$ (Conjecture \zcref{generalization_of_main_thm}).
For ${\bm x} \in P^{\vee}$, denote by $(\xi_{\bnu})_{\bm x}$ the restriction of $\xi_{\bnu}$ to $\cM \times \{ {\bm x} \}$. 
From the construction of $\xi_{\bnu}$, Theorem \zcref{main} and its proof imply the following theorem.
\begin{Thm}[Theorem \zcref{generic_orthogonal_property}]
Suppose $n=5$, and ${\bm  x}, {\bm y} \in P^{\vee} \setminus (\cup_{i = 1}^5 Z_i^{\pm} \cup \Delta)$. Then
\begin{equation*}
 H^i(\cM, (\xi_{\bnu})_{\bm x} \otimes (\xi_{\bnu})_{\bm y}) = 0
\end{equation*}
for any ${\bm x} \neq {\bm y}$, $i \geq 0$.\qed
\end{Thm}

We will also explain in the $ n=5$ case the computation of the cohomology of the structure sheaf of $\cM$ given by D. Arinkin \cite{AF} based on his discussions with R. Fedorov.
\begin{Thm}[D. Arinkin \cite{AF}, Corollary \zcref{proof_of_main2}]\label{main2}
Suppose $n = 5$. 
Under the hypothesis of Proposition \zcref{cohomology_of_bigraded_algebra},
 \begin{equation*}
  H^i(\cM, \cO_{\cM}) = \begin{cases}
                                     \bC & \text{if}\ i = 0,\\
                                     0 & \text{if}\ i > 0.
                                    \end{cases}
 \end{equation*}
\end{Thm}
In Proposition \zcref{delta_isomorphism}, we outline the proof of this assumption.
These theorems support the orthogonality property mentioned above. 
In future work, we will check the orthogonality property along the remaining locus, and prove the categorical equivalence (Conjecture \zcref{GLC-like}) in the $n = 5$ case.
Theorem \zcref{main2} is partially obtained by the author in \cite{M21a, M21b}, in which the statement is equivalent to a certain connecting map being an isomorphism.
Note that the methods given by Arinkin \cite{AF} work for arbitrary $n \geq 5$ if we assume the statements corresponding to Proposition \zcref{cohomology_of_compactification} and Lemma \zcref{canonical_bundle_over_M_H} in this paper.
We will explain the proof of Theorem \zcref{main} in \zcref{proof_of_main_theorem}.
The proof of Theorem \zcref{main2} is given in \zcref{Proof_of_Theorem2} as Corollary \zcref{proof_of_main2}.
\subsection{Proof of Theorem \zcref{main}}\label{proof_of_main_theorem}
To compute the cohomology of vector bundles on $\cM$, we construct its compactification. 
Deligne introduced a notion of $\lambda$-connections, and Simpson constructed a compactification of the moduli space of connections by using it (\cite{Sim91, Sim97}).
In \cite{A01}, Arinkin defined $\epsilon$-connections, a variant of Deligne--Simpson's $\lambda$-connections, and constructed a natural compactification of $\cM$.
While $\lambda$-connections give us an $\mathbb{A}^1$-family of moduli spaces, $\epsilon$-connections give us an $\mathbb{A}^1/\mathbb{G}_m$-family of them as explained below.

Suppose $E$ is a one-dimensional vector space, $\epsilon \in E$, $L$ is a rank $2$ vector bundle on $\bP^1$,
$\nabla \colon L \ra L \otimes \Omega^1_{\bP^1}(D) \otimes E$ is a $\bC$-linear map, and $\varphi \colon \bigwedge^2 L \xrightarrow{\sim} \cO_{\bP^1}$.

\begin{Def}\label{epsilon}
A collection $(L, \nabla, \varphi ; \epsilon \in E)$ is called an {\rm\bf{$\epsilon$-connection}} if the following conditions hold:
 \begin{itemize}
  \item[(1)] $\nabla(fs) = f \nabla s + s \otimes df \otimes \epsilon$ for $f \in \cO_{\bP^1}, s \in L$, 
  \item[(2)] $\varphi (\nabla s_1 \wedge s_2) + \varphi(s_1 \wedge \nabla s_2) = d(\varphi(s_1 \wedge s_2)) \otimes \epsilon$ for $s_1, s_2 \in L$,
  \item[(3)] The map $\res_{t_i}(\nabla) : L_{t_i} \ra (L \otimes \Omega^1_{\bP^1}(D) \otimes E)_{t_i} = L_{t_i} \otimes E$ induced by $\nabla$ has eigenvalues $\{ \epsilon \nu_i, - \epsilon \nu_i  \}$ for each $i$ $(1 \leq i \leq n)$,
  \item[(4)] $(L, \nabla)$ is irreducible; that is, there is no rank one subbundle $L_0 \subset L$ such that 
            $\nabla(L_0) \subset L_0 \otimes \Omega^1_{\bP^1}(D) \otimes E$.
 \end{itemize}
\end{Def}

Let $\overline{\cM}$ be the moduli stack of $\epsilon$-connections. 
Vector spaces $E$ for $\epsilon$-connections $(L, \nabla, \varphi; \epsilon \in E)$ form an invertible sheaf $\cE$ on $\overline{\cM}$ together with a natural section $\epsilon \in H^0(\overline{\cM}, \cE)$. 
Denote by $\cM_H \subset \overline{\cM}$ the closed substack defined by the equation $\epsilon = 0$.
Taking $E = \bC, \epsilon = 1$, we see that $\bnu$-$\mathfrak{s}l_2$-parabolic connections are particular cases of $\epsilon$-connections.
It is well-known that such connections are irreducible (cf. \cite[Proposition 1]{AL97}).
Moreover, if $\epsilon \neq 0$, there is a unique isomorphism $E \ra \bC$ such that $\epsilon \mapsto 1$.
It follows that the open substack $\overline{\cM} \setminus \cM_H$ corresponding to $\epsilon$-connections with $\epsilon \neq 0$ parametrizes all $\bnu$-$\mathfrak{s}l_2$-parabolic connections, and so, it is $\cM$.
Therefore, we have the map
\begin{equation*}
r: \overline{\cM} \ra \mathbb{A}^1/\mathbb{G}_m\ : \ (L, \nabla, \varphi; \epsilon \in E ) \mapsto [(\epsilon \in E)],
\end{equation*}
where the quotient stack $\mathbb{A}^1/\mathbb{G}_m$ is the moduli stack of pairs $(\epsilon \in E)$, and $r^{-1}([(0 \in E)]) = \cM_H$ and $r^{-1}([(\epsilon \in E)]) = \cM$ for $\epsilon \neq 0$. 

To show Theorem \zcref{main}, we need the next two propositions:
For $x \in \bP^1$, we also define the bundle $\xi_x$ over $\overline{\cM}$ whose fiber at $(L, \nabla, \varphi; \epsilon \in E)$ is $L_x$.
\begin{Prop}[Proposition \zcref{computation_along_boundary}]\label{pre_compactified_locus}
Suppose $x_1, \dots, x_4 \in \bP^1$ and $x_i \neq x_j$ for $i \neq j$. Then 
\begin{equation*}
 H^i(\cM_H, \xi_{x_1}\otimes \xi_{x_2} \otimes \xi_{x_3} \otimes \xi_{x_4} \otimes (\cE | _{\cM_H})^{\otimes k}) = 0,
\end{equation*}
for any $i, k$.
\end{Prop} 

\begin{Prop}[Proposition \zcref{computation_compactified}]\label{pre_compactified_vanish}
Suppose $x_1, \dots, x_4 \in \bP^1$. Then
\begin{equation*}
 H^i(\overline{\cM}, \xi_{x_1} \otimes \xi_{x_2} \otimes \xi_{x_3} \otimes \xi_{x_4}(-\cM_H)) = 0, 
 \end{equation*}
for any $i$.
\end{Prop}

We will show Proposition \zcref{pre_compactified_locus} in \zcref{proof_of_prop_1} and Proposition \zcref{pre_compactified_vanish} in \zcref{proof_of_prop_2} by decomposing vector bundles into line bundles.
For the $n = 4$ case, the space $M$, the underlying coarse moduli space of $\cM$, has long been studied as the Okamoto initial condition space of Painlev\'{e} V\hspace{-1.2pt}I equation \cite{Oka79}.
In particular, since $M$ is then an algebraic surface, we can construct and study its compactifications by using elementary algebro-geometric methods.
By contrast, when $n = 5$, the dimension of $M$ jumps to four, and this higher-dimensional version has received little attention.
Meanwhile, recent advances in the theory of apparent singular points enable us to introduce nice coordinates on $M$ (cf. \cite{LS}).
In this paper, we exploit these developments to provide an explicit geometric description of $\overline{\cM}$ given by the blowing-up of $\bP^2 \times (\bP^2)^{\vee}$ in the $n = 5$ case (Theorem \zcref{extended_app_bun}).
On the other hand, the boundary locus $\cM_H$ can be related to the moduli space of parabolic Higgs bundles.
By the general theory of Hitchin integrable systems, this space is isomorphic to a certain family of Jacobian varieties of spectral curves.
By using Arinkin's Fourier--Mukai transforms for compactified Jacobians \cite{A11, A13}, extended recently by Maulik--Shen--Yin in \cite{MSY} to the twisted case, we calculate the cohomology of vector bundles.

Denote by $j : \cM \hookrightarrow \overline{\cM}$ and $i : \cM_H \hookrightarrow \overline{\cM}$ the natural embeddings.
For a vector bundle $\cF$ on $\overline{\cM}$, we consider the filtration
\begin{equation*}
 \cF_0 := \cF \subset \cdots \subset \cF_k := \cF(k\cM_H) \subset \cdots \subset \cF_{\infty} :=j_*j^*\cF.
\end{equation*}
This yields $H^{\bullet}(\cM, \cF|_{\cM}) = H^{\bullet}(\overline{\cM}, \cF_{\infty}) = \varinjlim H^{\bullet}(\overline{\cM}, \cF_k)$.
Besides,
\begin{equation*}
\cF_k/\cF_{k-1} = i_*(\cF_k|_{\cM_H}) = i_*(\cF|_{\cM_H}\otimes (N_{\cM_H})^{\otimes k}),
\end{equation*}
where $N_{\cM_H} \simeq \cE|_{\cM_H}$ is the normal bundle to $\cM_H \subset \overline{\cM}$.
Therefore, we get the following lemma.

\begin{Lem}\label{isom_lemma}
Suppose $\cF$ be a vector bundle on $\overline{\cM}$ such that
\begin{equation*}
 H^{\bullet}(\cM_H, \cF|_{\cM_H} \otimes (N_{\cM_H})^{\otimes k}) = 0
\end{equation*}
for any $k > 0$.
Then, the natural maps $H^{\bullet}(\overline{\cM}, \cF) \ra H^{\bullet}(\cM, \cF|_{\cM})$ are isomorphisms.
\end{Lem}

\begin{proof}[Proof of Theorem \zcref{main}]
Set $\cF : = \xi_{x_1} \otimes \xi_{x_2} \otimes \xi_{x_3} \otimes \xi_{x_4}(-\cM_H)$.
Using Proposition \zcref{pre_compactified_locus} and Lemma \zcref{isom_lemma}, we get $H^{\bullet}(\overline{\cM}, \cF) = H^{\bullet}(\cM, \cF|_{\cM})$.
Now Proposition \zcref{pre_compactified_vanish} completes the proof.
\end{proof}

\subsection{Outline of the paper}
We briefly outline the contents of this paper.
In \zcref{chapter_jacobian_line_bundle}, we will compute the cohomology of torsion-free sheaves on twisted compactified Jacobians, which is needed to show  Proposition \zcref{pre_compactified_locus}.
In \zcref{M_H} and \zcref{M_bar}, we will show Proposition \zcref{pre_compactified_locus} and \zcref{pre_compactified_vanish} by studying the behavior of $\xi_x$ on $\cM_H$ and $\overline{\cM}$ respectively.
In \zcref{Proof_of_Theorem2}, we will show Theorem \zcref{main2} following the strategy given by Arinkin \cite{AF}.
In \zcref{Section_GLC-like}, we will explain the relationship between our main results and the geometric Langlands correspondence.

\subsection{Acknowledgements}
I would like to thank my Ph.D. supervisor, Vivek Shende, for his constant attention to this work and for his warm encouragement.
I am especially grateful to Dima Arinkin for sharing his idea and for the helpful discussions in the early stages of the project.
I am also very grateful to Michi-aki Inaba for his hospitality at Nara Women's University.
I am also very grateful to Frank Loray for his hospitality at Universit\'{e} de Rennes 1.
I also thank Tony Pantev, Junliang Shen, and Yukinobu Toda for many comments and suggestions.
The work presented in this article is supported by VILLUM FONDEN, VILLUM Investigator grant 37814, Novo Nordisk Foundation grant NNF20OC0066298, and academist crowdfunding.

\section{Cohomology of line bundles on relative twisted compactified Jacobians}\label{chapter_jacobian_line_bundle}
In this section, we calculate the cohomology of line bundles on relative twisted compactified Jacobians.
\subsubsection*{Notations}
We denote by $\cD_{qc}^b(X)$ (respectively, $\cD^b(X))$ the bounded derived category of quasi-coherent sheaves (respectively, coherent sheaves) over $X$. 

\subsection{Relative compactified Jacobians}
We recall here the Fourier-Mukai theorem of relative compactified Jacobians developed by Arinkin \cite{A11}.

Fix $g \geq 0$.
Let $p_C \colon C \rightarrow S$ be a family of projective integral curves with planar singularities of arithmetic genus $g$ over a base scheme $S$.
Let $J_C$ be the moduli space of pairs $(s, L)$, where $s \in S$ and $L$ is a degree $0$ line bundle on $C_s$.
Similarly, let $\overline{J}^0_C$ be the moduli space of pairs $(s, F)$, where $s \in S$ and $F$ is a degree $0$ torsion-free sheaf of generic rank one on $C_s$.
Then, $J_C \subset \overline{J}^0_C$ is an open subvariety.

Consider the Poincar\'e bundle $P^0$ on $J_C \times \overline{J}^0_C$.
Its fiber over $(s, L, F) \in J_C \times \overline{J}^0_C$ is equal to
\begin{equation}\label{fiber_formula1}
\det R\Gamma(C_s, F \otimes L) \otimes \det R\Gamma (C_s, \mathcal{O}) \otimes \det R\Gamma (C_s, F)^{-1} \otimes \det R\Gamma (C_s, L)^{-1}.
\end{equation}
More precisely, we can write $L \simeq \mathcal{O}_{C_s}(\sum a_i x_i)$ for a divisor $\sum a_i x_i$ supported by the smooth locus of $C_s$ and then 
\begin{equation}\label{fiber1}
P^0_{(s, L, F)} = \bigotimes (F_{x_i})^{\otimes a_i}.
\end{equation}
We can explain this normalization by using the universal line bundle $\cL$ (resp. the universal sheaf $\cF$) over $C \times _S J_C$ (resp. over $C \times_S \overline{J}^0_C$) and write
\begin{equation}\label{Arinkin_normalization}
 P^0 = P^0(\cL, \cF)
\end{equation}
to indicate the dependence of $P^0$ on $\cL$ and $\cF$.

Denote by $\mathfrak{j}$ the rank $g$ vector bundle on $S$ whose fiber over $s \in S$ is $H^1(C_s, \mathcal{O}_{C_s})$.
The relative dualizing sheaf for $q: J_C \rightarrow S$ equals $\Omega^g_{J/S} = q^*(\det (\mathfrak{j})^{-1})$.

\begin{Thm}[{\cite[Theorem5.1]{A11}}]\label{orthogonality}
Let $\pi_1 \colon J_C \times_S \overline{J}^0_C \rightarrow J_C$ be the projection. Then
\begin{equation*}
 \begin{split}
  R\pi_{1, *}P^0 &= (\Omega^g_{J/S})^{-1} \otimes \zeta_{0, *} \mathcal{O}_S[-g]\\
                       &=\zeta_{0, *} \det (\mathfrak{j}) [-g],
 \end{split}
\end{equation*}
where $\zeta_0 \colon S \rightarrow J_C$ is the zero section. \qed
\end{Thm}
In the proof of this theorem, Arinkin showed the next lemma.

\begin{Lem}\label{supp}
\begin{equation*}
 \Supp(R^{\bullet}\pi_{1, *} P^0) = \zeta_0(S).
\end{equation*} 
\end{Lem}
\begin{proof}
See the proof of \cite[Theorem5.1]{A11}.
\end{proof}
As a set, $\Supp(R^i \pi_{1, *}P^0)$ consists of pairs $(s, L) \in J_C$ such that the line bundle $L$ on $C_s$ satisfies $H^i((p^{0})^{-1}(s), P^0_L) \neq 0$. 
Here, $P^0_L$ is the restriction of $P^0$ to $\{(s, L)\} \times_S \overline{J}^0_C$ and $p^0 : \overline{J}^0_C \rightarrow S$.
Therefore, we get:

\begin{Thm}[{\cite[Theorem 1.2 (i)]{A11}}]\label{vanishing_degree_0}
Let $\zeta \colon U \rightarrow J_C$ be a local section over an open subset $U \subset S$.
If $\zeta(U) \cap \zeta_0(U) = \emptyset$, then we have $H^i(\overline{J}^0_{C/U}, P^0_{\zeta}) = 0$ for any $i$.
Here, $\overline{J}^0_{C/U}$ is the restriction of $\overline{J}^0_{C}$ to $(p^0)^{-1}(U)$ and
$P^0_{\zeta}$ is the restriction of $P^0$ to $\zeta(U) \times_S \overline{J}^0_C$. \qed
\end{Thm}

To understand this theorem, we will see the specific case.
Fix $s \in S$, and denote $\overline{\Pic}^0(C_s) := (p^0)^{-1}(s)$.
For a smooth point $x \in C_s$, let $\xi_x$ be the line bundle on $\overline{\Pic}^0(C_s)$ whose fiber over $F$ equals $F_x$.
By using \zcref{fiber1}, we get the following corollary.

\begin{Cor}\label{cohomology_degree0Jac}
If $L \simeq \cO_{C_s}(\sum a_i x_i) \not\simeq \cO_{C_s}$, then
\begin{equation*}
 H^k(\overline{\Pic}^0(C_s), \bigotimes (\xi_{x_i})^{\otimes a_i}) = 0
\end{equation*}
for any $k$. \qed
\end{Cor}

Theorem \zcref{orthogonality} can be formulated in terms of the Fourier-Mukai functor
\begin{equation*}
\mathfrak{F} \colon D^b_{qc}(J_C) \rightarrow D^b_{qc}(\overline{J}^0_C)\ : \ \mathcal{F} \mapsto R\pi_{2, *}(\pi_1^*(\mathcal{F}) \otimes P^0)
\end{equation*}
given by $P^0$.

\begin{Thm}[{\cite[Theorem 1.4 (ii)]{A11}}]\label{fullyfaithful}
$\mathfrak{F}$ is fully-faithful. \qed
\end{Thm}

By theorem \zcref{fullyfaithful} and Koszul complex, we can compute the cohomology of structure sheaf $\mathcal{O}_{\overline{J}^0_C}$.

\begin{Cor}[{\cite[Theorem 1.2 (ii), Proposition 6.1]{A11}}]
We have an isomorphism of graded algebras
\begin{equation*}
R^{\bullet} p^0_* \mathcal{O}_{\overline{J}^0_C} = \bigwedge^{\bullet} R^1p_{C, *} \mathcal{O}_C.
\end{equation*}
Here, $p^0 \colon \overline{J}^0_C \rightarrow S$. \qed
\end{Cor}

\subsection{Autoduality of compactified Jacobians}

We summarize here the results of Arinkin \cite{A13} that extend the Theorem \zcref{fullyfaithful} as an autoequivalence of $ D^b_{qc}(\overline{J}^0_C)$.

In \cite{A13}, Arinkin constructed the Poincar\'e sheaf $\overline{P}^0$ over $\overline{J}^0_C \times_S \overline{J}^0_C$.
Let $j \colon J_C \times \overline{J}^0_C \cup \overline{J}^0_C \times J_C \hookrightarrow \overline{J}^0_C \times \overline{J}^0_C$ be an open embedding and $\overline{\pi}_i \colon \overline{J}^0_C \times \overline{J}^0_C \ra \overline{J}^0_C$, $(i = 1, 2)$ be the projection.

\begin{Thm}[{\cite[Theorem A, Lemma 6.1]{A13}}]\label{Arinkin_existence_thm}
There exists a coherent sheaf $\overline{P}^0$ on $\overline{J}^0_C \times_S \overline{J}^0_C$ with the following properties:
\begin{itemize}
 \item[(1)]$\overline{P}^0 = j_*P^0$,
 \item[(2)]$\overline{P}^0$ is flat for the projection $\overline{\pi}_2 \colon \overline{J}^0_C \times_S \overline{J}^0_C \ra \overline{J}^0_C$,
 \item[(3)]$\overline{P}^0$ is a maximal Cohen-Macaulay sheaf on $\overline{J}^0_C \times_S \overline{J}^0_C$.
\end{itemize} \qed
\end{Thm}

The Poincar\'e sheaf $\overline{P}^0$ provides a categorical autoduality of $\overline{J}^0_C$.
\begin{Thm}[{\cite[Theorem C]{A13}}]\label{arinkin_autodual}
The Fourier-Mukai functor
\begin{equation}
\overline{\mathfrak{F}} \colon D^b_{qc}(\overline{J}^0_C) \ra D^b_{qc}(\overline{J}^0_C) \ : \ \cG \mapsto R\overline{\pi}_{1,*}(\overline{\pi}_2^*(\cG)\otimes \overline{P}^0)
\end{equation}
is an equivalence of categories. \qed
\end{Thm}

Theorem \zcref{arinkin_autodual} comes from the next proposition:
Set 
\begin{equation*}
\Psi := Rp_{13,*}(p_{12}^* (\overline{P}^0)^{\vee} \otimes p_{23}^* \overline{P}^0) \in D^b_{qc}(\overline{J}^0_C \times \overline{J}^0_C).
\end{equation*}
Here $(\overline{P}^0)^{\vee} := {\mathcal{H}om}(\overline{P}^0, \cO_{\overline{J}^0_C \times \overline{J}^0_C})$.
Denote the projection $\overline{J}^0_C \times_S \overline{J}^0_C \ra S$ by $\pi$ and diagonal in $\overline{J}^0_C \times \overline{J}^0_C$ by $\Delta$.
Recall that $\mathfrak{j}$ is the rank $g$ vector bundle on $S$ whose fiber over $s \in S$ is $H^1(C_s, \mathcal{O}_{C_s})$.

\begin{Prop}[{\cite[Proposition 7.1]{A13}}]
\begin{equation*}
 \Psi \simeq \cO_{\Delta}[-g] \otimes \pi^* \det(\mathfrak{j}).
\end{equation*} \qed
\end{Prop}

This proposition also implies the next statement, which is similar to Theorem \zcref{orthogonality}.
\begin{Prop}\label{support_of_poincare_sheaf}
\begin{equation*}
  R\overline{\pi}_{1, *}\overline{P}^0 = \overline{\zeta}_{0, *} \det (\mathfrak{j}) [-g],
\end{equation*}
where $\overline{\zeta}_0 \colon S \rightarrow \overline{J}^0_C$ is the zero section.
Especially, $\Supp(R\overline{\pi}_{1, *}\overline{P}^0) = \overline{\zeta}_0(S)$. \qed
\end{Prop}

\subsection{Relative twisted compactified Jacobians}
In \cite{MSY}, Maulik, Shen, and Yin extended Arinkin's results to the twisted case.
We assume that the total space of $C \ra S$ is nonsingular, and there is a multi-section $\sigma \colon S \ra C$, $D := \sigma(S)$ of degree $r$ which is finite and flat over $S$.

Let $\overline{J}^d_C$ be the moduli space of pairs $(s, F)$, where $s \in S$ and $F$ is a degree $d$ torsion-free sheaf of generic rank one on $C_s$.
We assume that for any degree $d$, the relative twisted compactified Jacobian $\overline{J}^d_C$ is a nonsingular quasi-projective variety.
We denote the natural projection map by
\begin{equation*}
 p^d \colon \overline{J}^d_C \ra S.
\end{equation*}

\subsubsection{Trivialization along a multi-section}\label{m2_gerbe}
For any $S$-scheme $T$, we consider a flat family $\cF^d_T$ over $C \times_S T$ of rank $1$ torsion-free sheaves of degree $d$ on the curves parametrized by $T$.
We define 
\begin{equation*}
 \mathcal{R}^d_T := \det (p_{T, *}(\cF^d_T|_{D \times_S T})) \in \Pic (T),
\end{equation*}
where $p_T \colon D \times_S T \ra T$ is the natural projection.
We say that this family over $T$ is trivialized along the multi-section $D \subset C$, if there is a specified isomorphism
\begin{equation*}
 \mathcal{R}^d_T \simeq \cO_T \in \Pic (T).
\end{equation*}

Let $\overline{\cJ}^d_C$ be the functor sending any $S$-scheme $T$ to a groupoid given by the data
\begin{equation*}
 \cF^d_T \leadsto C \times_S T 
\end{equation*}
satisfying the same conditions as for the stack of the degree $d$ compactified Jacobian, with an extra assumption that $\cF^d_T$ is trivialized along the multi-section $D$.

\begin{Prop}[{\cite[Proposition 4.1]{MSY}}]
 The functor $\overline{\cJ}^d_C$ is represented by a Deligne-Mumford stack which is a $\mu_r$-gerbe over $\overline{J}^d_C$.\qed
\end{Prop}

We have got, for any $d$, a nonsingular Deligne-Mumford stack $\overline{\cJ}^d_C$ which is a $\mu_r$-gerbe over $\overline{J}^d_C$, together with the universal family $\cF^d$ of rank $1$ degree $d$ torsion-free sheaves on $C \times_S \overline{\cJ}^d_C$, trivialized along the multi-section $D$;
\begin{equation*}
 \det (p_{\cJ, *}(\cF^d|_{D \times_S \overline{\cJ}^d_C})) \simeq \cO_{\overline{\cJ}^d_C} \in \Pic (\overline{\cJ}^d_C),
\end{equation*}
where $p_{\cJ} \colon D \times_S \overline{\cJ}^d_C \ra \overline{\cJ}^d_C$ is the natural projection.
 
Let $\Coh(\overline{\cJ}^d_C)_{(k)}$ (resp. $D^b(\overline{\cJ}^d_C)_{(k)}$) be the full subcategory of $\Coh(\overline{\cJ}^d_C)$ (resp. $D^b(\overline{\cJ}^d_C)$) consisting of objects for which the action of $\mu_r$ on fibers is given by the character $\lambda \mapsto \lambda^k$ of $\mu_r$.
Let us consider two integers $d, k$.
We put $\cF^k, \cF^d$ in Arinkin's formula \zcref{Arinkin_normalization} and obtain a Poincar\'e sheaf
\begin{equation*}
 \overline{P}^{(k, d)} := P^0(\cF^k, \cF^d) \in \Coh (\overline{\cJ}^k_C \times_S \overline{\cJ}^d_C)_{(d, k)}.
\end{equation*}
We also define
\begin{equation*}
 (\overline{P}^{(k, d)})^{-1} := {\mathcal{H}om}_{\overline{\cJ}^k_C \times_S \overline{\cJ}^d_C}(\overline{P}^{(k, d)}, p^*_2 \Omega^g_{\overline{\cJ}^d_S/S})[g] \in D^b(\overline{\cJ}^k_C \times_S \overline{\cJ}^d_C)_{(-k, -d)},
\end{equation*}
where $\Omega^g_{\overline{\cJ}^d_S/S}$ is the relative dualizing sheaf with respect to $p^d : \overline{\cJ}^d_C \ra S$.

Because $\overline{P}^{(k, d)}$ lies in the isotypic category, the Fourier-Mukai transform
\begin{equation*}
 \overline{\mathfrak{F}}^{(k, d)} : D^b(\overline{\cJ}^k_C) \ra D^b(\overline{\cJ}^d_C)
\end{equation*}
with the kernel $\overline{P}^{(k, d)}$ is only non-zero on the following isotypic components;
\begin{equation*}
 \overline{\mathfrak{F}}^{(k, d)} : D^b(\overline{\cJ}^k_C)_{(-d)} \ra D^b(\overline{\cJ}^d_C)_{(k)}.
\end{equation*}

\begin{Thm}[{\cite[Proposition 4.2]{MSY}}]\label{thm_general_Fourier_Mukai}
The Fourier-Mukai functor
\begin{equation}\label{general_Fourier_Mukai}
\overline{\mathfrak{F}}^{(k, d)} \colon D^b(\overline{\cJ}^k_C)_{(-d)} \ra D^b(\overline{\cJ}^d_C)_{(k)} \ : \ \cG \mapsto R\overline{\pi}_{1,*}(\overline{\pi}_2^*(\cG)\otimes \overline{P}^{(k, d)})
\end{equation}
is an equivalence of categories.
Its quasi-inverse is given as the Fourier-Mukai functor with the kernel $(\overline{P}^{(k, d)})^{-1}$. \qed
\end{Thm}

\subsubsection{\'Etale local descriptions}
Let $U$ be a \'etale neighborhood of the base $S$.
We may assume that $C \ra U$ admits simultaneously a section $\gamma$ and a multi-section $\sigma$;
\begin{equation*}
 U \simeq \gamma(U) \subset C \ra U, \ D:= \sigma(U) \subset C \ra U.
\end{equation*}
They are independent and do not have any non-trivial relation.
Due to the existence of the section, the relative compactified Jacobian $\overline{J}^d_C$ are identified for any choice of $d$.

We will compare $\cF^d$ on $C \times_U \overline{\cJ}^d_C$ which is trivialized along $D$ with the normalized universal sheaf $\cF$ on $C \times_U \overline{J}^0_C$ which is trivialized along the section $\gamma(U)$.

\begin{Prop}[{\cite[Proposition 4.3]{MSY}}]
 There is a $U$-morphism
 \begin{equation*}
  \iota_d : \overline{\cJ}^d_C \ra \overline{J}^0_C
 \end{equation*}
satisfying the following properties:
 \begin{itemize}
  \item[(1)] We have 
     \begin{equation*}
      \cF^d \simeq (\id_C \times_C \iota_d)^*\cF \otimes p_C^*\cO_C(d \gamma(U)) \otimes p_{\cJ}^*\cL_d,
     \end{equation*}
     where $\cL_d \in \Pic(\overline{\cJ}^d_C)$ and $p_C, p_{\cJ}$, are the natural projections from $C \times_U \overline{\cJ}^d_C$.
  \item[(2)]For any $s \in U$ with $C_s$ a non-singular curve, the restriction of $\cL_d$ to the fiber $\cJ^d_{C_s}$ has trivial first Chern class in $H^2(\cJ^d_C, \mathbb{Q})$.
 \end{itemize} \qed
\end{Prop}

We consider the morphism
\begin{equation*}
 \iota_k \times_U \iota_d : \overline{\cJ}^k_C \times_U \overline{\cJ}^d_C \ra \overline{J}^0_C \times_U \overline{J}^0_C.
\end{equation*}
The normalized Poincar\'e sheaf $\overline{P}^0$ and its inverse $(\overline{P}^0)^{-1}$ are canonically defined over $\overline{J}^0_C \times_U \overline{J}^0_C$ (Theorem \zcref{Arinkin_existence_thm}).

\begin{Thm}[{\cite[Corollary 4.4]{MSY}}]
 We have 
 \begin{equation}\label{bundle_relation}
  \overline{P}^{(k, d)} \simeq (\iota_k \times_U \iota_d)^*\overline{P}^0 \otimes (L_k^{\otimes d} \boxtimes L_d^{\otimes k}),
 \end{equation}
 \begin{equation}
  (\overline{P}^{(k, d)})^{-1} \simeq (\iota_k \times_U \iota_d)^*(\overline{P}^0)^{-1} \otimes (L_k^{\vee \otimes d} \boxtimes L_d^{\vee \otimes k}).
 \end{equation}
Here $L_d, L_k$ are line bundles over $\overline{\cJ}^d_C, \overline{\cJ}^k_C$ respectively.
Moreover, for any $s \in U$ with $C_s$ non-singular, the restrictions of those line bundles to the fibers over$s$ have homologically trivial first Chern classes. \qed
\end{Thm}

\subsection{Cohomology of line bundles on relative twisted compactified Jacobians}
Firstly, we consider the special case of the equivalence \zcref{general_Fourier_Mukai};
\begin{equation*}
\overline{\mathfrak{F}}^{(0, d)} \colon D^b(\overline{\cJ}^0_C)_{(-d)} \ra D^b(\overline{\cJ}^d_C)_{(0)} \ : \ \cG \mapsto R\overline{\pi}_{1,*}(\overline{\pi}_2^*(\cG)\otimes \overline{P}^{(0, d)}).
\end{equation*}
Here, $D^b(\overline{\cJ}^d_C)_{(0)}$ is equivalent to $D^b(\overline{J}^d_C)$. 
Denote by $\mathcal{N}^d$ the sheaf on $\overline{\cJ}^0_C$ that is the structure sheaf of the zero section equipped with the action of $\mu_r$ with weight $d$.
Because of the normalization of the Poincar\'e bundle, we have
\begin{equation*}
 \overline{\mathfrak{F}}^{(0, d)} (\mathcal{N}^{-d}) = \mathcal{O}_{\overline{J}^d_C}.
\end{equation*}
Therefore, theorem \zcref{thm_general_Fourier_Mukai} implies the following formula:

\begin{Cor}\label{direct image}
\begin{equation}
Rp^d_*\mathcal{O}_{\overline{J}^d_{C}} \simeq R\mathcal{H}om_S (\mathcal{N}^{-d}, \mathcal{N}^{-d}).
\end{equation}
Here the object on the right-hand side is the pushforward of $R\mathcal{H}om_{\overline{\cJ}^0_C}(\mathcal{N}^{-d}, \mathcal{N}^{-d})$ to $S$.
\end{Cor}
\begin{proof}
For the projections $p^d \colon \overline{J}^d_C \rightarrow S$ and $q \colon \overline{\cJ}^0_C \rightarrow S$, we have
 \begin{equation*}
  \begin{split}
  Rq_* R \mathcal{H}om _{\overline{\cJ}^0_C}(\mathcal{N}^{-d}, \mathcal{N}^{-d}) & \simeq Rp^d_* \overline{\mathfrak{F}}^{(0, d)}(R \mathcal{H}om_{\overline{\cJ}^0_C}(\mathcal{N}^{-d}, \mathcal{N}^{-d}))\\
  &\simeq Rp^d_* R\mathcal{H}om_{\overline{J}^d_C}(\overline{\mathfrak{F}}^{(0, d)}(\mathcal{N}^{-d}), \overline{\mathfrak{F}}^{(0, d)}(\mathcal{N}^{-d}))\\
  &\simeq Rp^d_* R \mathcal{H}om_{\overline{J}^d_C}(\mathcal{O}_{\overline{J}^d_C}, \mathcal{O}_{\overline{J}^d_C})\\
  &\simeq Rp^d_* \mathcal{O}_{\overline{J}^d_C}.
  \end{split}
 \end{equation*}
\end{proof}

The Koszul complex allows us to compute the cohomology of the right hand side of Corollary \zcref{direct image}.

\begin{Cor}\label{higher_direct_image_of_strucuture_sheaf}
We have an isomorphism of graded algebras
\begin{equation*}
 R^{\bullet}p^d_* \mathcal{O}_{\overline{J}^d_C} = \bigwedge^{\bullet} R^1 p_{C, *} \mathcal{O}_C.
\end{equation*}\qed
\end{Cor}

Next, for general integers $d, k$, Proposition \zcref{support_of_poincare_sheaf} and \zcref{bundle_relation} give us the next statement:
\begin{Lem}
 \begin{equation}
  \Supp(R\overline{\pi}_{k, *}\overline{P}^{(k, d)}) = \iota_k^{-1}(\overline{\zeta}_0(U)).
 \end{equation}
 Here $\overline{\pi}_k \colon \overline{\cJ}^k_C \times_U \overline{\cJ}^d_C \ra \overline{\cJ}^k_C$ is the natural projection. \qed
\end{Lem}
This lemma implies the next statement, which is the generalization of Theorem \zcref{vanishing_degree_0}.

\begin{Thm}
Let $\overline{\zeta} \colon U \rightarrow \overline{\mathcal{J}}^k_C$ be a local section over an \'etale neighborhood $U \subset S$.
If $\overline{\zeta}(U) \cap \iota_k^{-1}(\overline{\zeta}_0(U)) = \emptyset$, then we have $H^i(\overline{\mathcal{J}}^d_{C/U}, \overline{P}^{(k, d)}_{\overline{\zeta}}) = 0$ for any $i$.
Here, $\overline{\mathcal{J}}^d_{C/U}$ is the restriction of $\overline{\mathcal{J}}^d_{C}$ to $(p^d)^{-1}(U)$ and
$\overline{P}^{(k, d)}_{\overline{\zeta}}$ is the restriction of $\overline{P}^{(k, d)}$ to $\overline{\zeta}(U) \times_S \overline{\mathcal{J}}^d_C$. \qed
\end{Thm}

\begin{Cor}\label{jacobian_vanishing_of_torsion_free_sheaves}
For $F \in \overline{\Pic}^k(C_s)$ such that $F \not\simeq  \cO_{C_s}(k \gamma(s))$, we have 
\begin{equation}
H^i(\overline{\Pic}^d(C_s), \overline{P}^{(k, d)}_F) = 0,
\end{equation}
for any $i$.
Here, $\overline{P}^{(k, d)}_F$ is the restriction of $\overline{P}^{(k, d)}$ to $\{(s, F)\} \times_S \overline{\cJ}^d_C$. \qed
\end{Cor}

\begin{Rem}[{\cite[Remarks (iv)]{A11}}]
\rm{In \zcref{chapter_jacobian_line_bundle},
we assumed the integrality of $C$ to avoid working with stability conditions for sheaves on $C$.
If one fixes an ample line bundle on $C$ and defines the compactified Jacobian $\overline{J}^d_C$ to be the moduli space of semi-stable torsion-free sheaves of degree $d$ of generic rank one, our argument works well.}
\end{Rem}

\section{Bundles $\xi_x$ on $\mathcal{M}_H$}\label{M_H}
In this section, we describe the boundary locus $\cM_H$, and prove Proposition \zcref{pre_compactified_locus} as Proposition \zcref{computation_along_boundary} by using Corollary \zcref{jacobian_vanishing_of_torsion_free_sheaves}.

\subsection{Relation between $\cM_H$ and ${\mathfrak{s}l_2}$-Higgs bundles}
Let us fix a set of $n$-distinct points $\bm{t} = \{ t_1, \dots, t_n \}$ on the Riemann sphere $\bP^1$ and define the divisor $D := t_1 + \cdots + t_n$.
Here, we study the relation between $\cM_H$ and the moduli space of ${\mathfrak{s}l_2}$-Higgs bundles over $(\bP^1, D)$ that we will define below.

Let $(L, \nabla, \varphi ; \epsilon \in E) \in \overline{\cM}$ be an $\e$-connection, and $U$ be the formal disk centered at $t_i$.
Trivializing $L|_U$, we can write
\begin{equation*}
\nabla|_U = \e d + A, \ A \in \mathfrak{s}l_2 \otimes \Omega^1_U(t_i) \otimes_{\bC} E.
\end{equation*}
Then, $\tr (A)$ and $\det(A)$ are well-defined, i.e., independent of the trivialization, as sections of 
$(\Omega^1_{\bP^1}(t_i)/\Omega^1_{\bP^1}) \otimes_{\bC} E$ and $(\Omega^{\otimes 2}_{\bP^1}(2t_i)/\Omega^{\otimes 2}_{\bP^1}(t_i)) \otimes _{\bC} E^{\otimes 2}$ respectively.
Performing at every $t_i$, we get well-defined sections of $(\Omega^1_{\bP^1}(D)/\Omega^1_{\bP^1}) \otimes_{\bC} E$ and $(\Omega^{\otimes 2}_{\bP^1}(2D)/\Omega^{\otimes 2}_{\bP^1}(D)) \otimes _{\bC} E^{\otimes 2}$, which we denote $[ \tr(\nabla)]$ and $[ \det(\nabla) ]$ respectively.

 For $(L, \nabla, \varphi; 1 \in \bC) \in \cM$, in a suitable trivialization of $L$ over the formal disk $U$ at $t_i$, we can write $\nabla$ as
 \begin{equation*}
  \begin{split}
    \nabla|_U &= d + a \frac{dz}{z}, \ a \in \mathfrak{s}l_2\\
                   &= d + \begin{pmatrix}
                                 \omega_i^+ & 0 \\
                                 0 & \omega_i^-
                               \end{pmatrix},
  \end{split}
 \end{equation*}
where $\omega_i^+, \omega_i^-$ are $1$-forms on the formal disk.
Then, we have
\begin{equation*}
[\tr (\nabla)] = \lambda_1 := (\omega_i^+ + \omega_i^-), \ \ [\det(\nabla)] = \lambda_2 := (\omega_i^+ \omega_i^-).
\end{equation*}
 For a general $\e$-connection $(L, \nabla, \varphi; \epsilon \in E) \in \overline{\cM}$, we have
 \begin{equation*}
[\tr (\nabla)] = \lambda_1 \otimes \e, \ \ [\det(\nabla)] = \lambda_2 \otimes \e^{\otimes 2}.
\end{equation*}
In particular for $\e = 0$ case, $(L, \nabla, \varphi; 0 \in E) \in \cM_H$ satisfies 
\begin{equation*}
\tr (\nabla) \in H^0(\bP^1, \Omega^1_{\bP^1}) \otimes E = 0,  \ \det(\nabla) \in H^0(\bP^1, \Omega^{\otimes 2}(D)) \otimes E^{\otimes 2}.
\end{equation*}

\begin{Lem}[{\cite[p. 215 Example]{A01}}]\label{det_condition}
$(L, \nabla, \varphi; 0 \in \bC)$ is an $\epsilon$-connection if and only if $\det (\nabla) \neq 0$. 
\end{Lem}
\begin{proof}
If $\det(\nabla) = 0$, then any rank one subbundle $L_0 \subset \Ker(\nabla)$ is $\nabla$-invariant.
Therefore $(L, \nabla)$ is reducible.
Conversely, assume that $(L, \nabla)$ is reducible. 
Then $\nabla$ has eigenvalues $\omega_{\pm} \in H^0(\bP^1, \Omega^1_{\bP^1}(D))$.
We have $\omega_+ + \omega_- = \tr(\nabla) = 0$ and $\omega_+ \omega_- = \det(\nabla) \in H^0(\bP^1, \Omega^{\otimes 2}(D))$.
Hence $\omega_{\pm} \in H^0(\bP^1, \Omega^1_{\bP^1}) = 0$, this implies that $\det(\nabla) = \omega_+ \omega_- = 0$.
Clearly other conditions (1)-(3) of Definition \zcref{epsilon} are satisfied.
\end{proof}

Therefore, $\cM_H$ is the moduli stack of $(L, \nabla, \varphi; 0 \in E)$, where
\begin{itemize}
 \item[(1)] $L$ is a rank $2$ vector bundle on $\bP^1$,
 \item[(2)] $E$ is a one dimensional vector space,
 \item[(3)] $\varphi \colon \bigwedge^2 L \xrightarrow{\sim} \cO_{\bP^1}$, 
 \item[(4)] $\nabla \colon L \ra L \otimes \Omega^1_{\bP^1}(D) \otimes E$ is an $\cO_{\bP^1}$-linear homomorphism,
 \item[(5)] the map $\res_{t_i}(\nabla) : L_{t_i} \ra L_{t_i} \otimes E$ induced by $\nabla$ is a nilpotent matrix,
 \item[(6)] $\tr(\nabla) = 0$,
 \item[(7)] $\det(\nabla) \in H^0(\bP^1, \Omega^{\otimes 2}(D)) \otimes E^{\otimes 2}$, and $\det(\nabla) \neq 0$.
\end{itemize}

We say that a triplet $(L, \nabla, \varphi)$ is {\it{holomorphic}} at $t \in \{t_1, \dots, t_n \}$ if $\res_{t}(\nabla) = 0$.

Now, let us introduce the corresponding Higgs bundles.

\begin{Def}\label{sl2_Higgs}
A {\rm\bf{$\mathfrak{s}l_2$-Higgs bundle}} over $(\bP^1, D)$ is a collection $(L, \theta, \varphi)$, where
\begin{itemize}
 \item[(1)] $L$ is a rank $2$ vector bundle on $\bP^1$,
 \item[(2)] $\varphi \colon \bigwedge^2 L \xrightarrow{\sim} \cO_{\bP^1}$, 
 \item[(3)] $\theta \colon L \ra L \otimes \Omega^1_{\bP^1}(D)$ is a traceless $\cO_{\bP^1}$-linear homomorphism,
 \item[(4)] the map $\res_{t_i}(\theta) \in \End (L_{t_i})$ induced by $\theta$ is a nilpotent matrix.
\end{itemize}
\end{Def}
We call $\theta$ {\it{Higgs field}}.
Let $\cY$ be the moduli stack of $\mathfrak{s}l_2$-Higgs bundles over $(\bP^1, D)$ with the condition $\det(\theta) \in H^0(\bP^1, \Omega^{\otimes 2}(D))$,  $\det(\theta) \neq 0$.
Note that, similarly to Lemma \zcref{det_condition}, every $\mathfrak{s}l_2$-Higgs bundle $(L, \theta, \varphi)$ over $(\bP^1, D)$ with $\det(\theta) \neq 0$ is irreducible, that is, there is no rank one subbundle $L_0 \subset L$ such that $\theta(L_0) \subset L_0 \otimes \Omega^1_{\bP^1}(D)$ (see also {\cite[Proposition 3.1 and 3.2]{FL23a}}).

The multiplicative group $\mathbb{G}_m$ acts on $\cY$ via
\begin{equation*}
a \cdot (L, \theta, \varphi) = (L, a\theta, \varphi), \ (L, \theta, \varphi) \in \cY, a \in \mathbb{G}_m.
\end{equation*}
Let us fix an isomorphism $\mu \colon E \xrightarrow{\sim} \bC$ which induces 
\begin{equation*}
 \id \otimes \mu^{\otimes 2} \colon H^0(\bP^1, \Omega^{\otimes 2}_{\bP^1}(D)) \otimes E^{\otimes 2} \xrightarrow{\sim} H^0(\bP^1, \Omega^{\otimes 2}_{\bP^1}(D)).
\end{equation*}
Then, we can identify $\cM_H$ with the quotient stack $\cY/\mathbb{G}_m$.
Consider the corresponding quotient map 
\begin{equation*}
\pi \colon \cY \ra \cM_H \ : \ (L, \theta, \varphi) \mapsto (L, \nabla, \varphi; 0 \in \bC).
\end{equation*}
We can easily check that
$\pi^* (\cE | _{\cM_H}) \simeq \cO_{\cY}$ (see \zcref{section_cohomology_of_hitchin_systems}). 

\subsection{Spectral curves and BNR spectral correspondence}
Put $V := H^0(\bP^1, \Omega^{\otimes 2}_{\bP^1}(D))$ and $S := V \setminus \{0\}$.
For any $s \in S$, we consider the spectral curve $C_s$ in the total space of $\Omega^1_{\bP^1}(D)$ as follows (cf. \cite[Section 3]{BNR}, \cite[Section 2.3]{FL23a}):
We can define a structure of commutative ring on $\cO_{\bP^1} \oplus (\Omega^1_{\bP^1}(D))^{-1}$ induced by $s$:
\begin{equation*}
(f_1, \omega_1) \times (f_2, \omega_2) := (f_1f_2 - s \otimes \omega_1 \otimes \omega_2, f_1 \omega_2 + f_2 \omega_1),
\end{equation*}
and this makes it an $\cO_{\bP^1}$-algebra.
It will be denoted by $\cA_s$, and is locally given by 
\begin{equation*}
\cA_s(U) = \frac{\cO_{\bP^1}(U)[t]}{(t^2 + s)}.
\end{equation*}
Then the spectral curve is given by $C_s := \Spec (\cA_s)$.
Note that for any $s$, we have a natural projection $\pi_s \colon C_s \ra \bP^1$ as a degree two map, and the push forward of $\cO_{C_s}$ to $\bP^1$ is identified with $\cO_{\bP^1} \oplus (\Omega^1_{\bP^1}(D))^{-1}$. 

\begin{Lem}\label{C_s_integral}
$C_s$ is an integral curve for any $s \in S$.  
\end{Lem}
\begin{proof}
Suppose that $C_s$ is non-integral for $s \in S$.
Then,
$P_s(t) = t^2 + s$,
the corresponding characteristic polynomial of $\theta$, is reducible over the function field of $\bP^1$.
Hence, on a Zariski open subset $U_i$ containing $t_i$, we can write
\begin{equation*}
P_s(t) = (t - a_i)(t - b_i)
\end{equation*} 
with $a_i, b_i \in \Omega^1_{\bP^1}(D)(U_i)$.
Then, it is easy to check that these $a_i$ glue together and give a global section $a \in H^0(\bP^1, \Omega^1_{\bP^1}(D))$ such that $s = -a^2$.

We now use the fact that $s \in H^0(\bP^1, \Omega_{\bP^1}^{\otimes 2}(D)) \setminus \{0 \}$.
If $a$ had a pole at one of the marked points $t_i$, then $a^2$ would have a pole of order at least two at $t_i$.
This is impossible, because $s$ has at most a simple pole along $D$.
Therefore, $a$ has no pole along $D$.
This means that $a \in H^0(\bP^1, \Omega^1_{\bP^1}) = 0$.
Thus $a = 0$, and hence $s = 0$, contradicting the assumption that $s \neq 0$.
Therefore, $C_s$ is integral.
\end{proof}

Denote by $\sigma_s \colon C_s \ra C_s$ the involution induced by
\begin{equation*}
 \sigma_s^* \colon \cA_s \ra \cA_s \ :\ (f, \omega) \mapsto (f, -\omega).
\end{equation*}
For an invertible sheaf $F$ on $C_s$, there is a natural action of $\sigma_s$ on the sheaf $F \otimes \sigma_s^*F$.
So, there is a natural invertible sheaf $\norm(F)$ on $\bP^1$ such that $F \otimes \sigma_s^*F = \pi_s^*\norm(F)$. 
Moreover, $\bigwedge^2(\pi_{s, *}F) = \bigwedge^2(\pi_{s, *}\cO_{C_s}) \otimes \norm(F) = \bigwedge^2 \cA_s \otimes \norm(F) = (\Omega^1_{\bP^1}(D))^{-1} \otimes \norm(F)$.
Here, $\norm$ is the norm map (cf. \cite[Chapter I\hspace{-1.2pt}V, Ex.  2.6]{Ha}).
This computation can be extended to torsion free rank one sheaves on $C_s$.

Let $\overline{\mathcal{J}}^{n-2}_C$ be the moduli stack of pairs $(s, F)$, where $s \in S$ and $F$ is a degree $n-2$ torsion-free sheaf of generic rank one on $C_s$.
There is a projection map $p \colon \overline{\mathcal{J}}^{n-2}_C \ra S$ and we denote $\overline{\Pic}^{n-2}(C_s) := p^{-1}(s)$.
There is a relation between $\cY$ and $\overline{\cJ}^{n-2}_C$, the so-called BNR spectral correspondence (cf. \cite[Section 3]{BNR}, \cite[Section 2.3]{FL23a}).

\begin{Prop}\label{BNR}
$\cY$ is naturally isomorphic to $\overline{\cJ}^{n-2}_C$. 
\end{Prop}
\begin{proof}
From Lemma \zcref{C_s_integral}, $C_s$ is integral for all $s \in S$.
Let $(L, \theta, \varphi)$ be a point of $\cY$.
Then $s = \det(\theta) \in S$ and $L$ is an $\cA_s$-module with respect to the multiplication
\begin{equation*}
(f, \omega)e := fe + \omega \otimes \theta(e)
\end{equation*}
for $e \in L$ and $(f, \omega) \in \cA_s$. It defines a sheaf $F$ of $\cA_s$-module on $C_s$.
Since $C_s$ is integral, then $F$ is torsion-free.

The inverse construction is given by $F \mapsto L := (\pi_s)_*F$.
Here $L$ is equipped with an $\cO_{\bP^1}$-linear map
\begin{equation*}
 \Theta \colon \cA_s \times L \ra L.
\end{equation*}
The section $t$ of $\pi_s^*(\Omega_{\bP^1}(D))$ gives a family of endomorphisms
\begin{equation*}
 \theta_U \colon L(U) \ra L(U) \ : \ e \mapsto \Theta(t, e) 
\end{equation*}
for each open set $U \subset \bP^1$ and which glue together to give a global homomorphism $\theta \colon L \ra L \otimes \Omega_{\bP^1}(D)$ satisfying $\theta^2 + \det(\theta)I_2 = 0$.
As explained above, $\varphi$ induces $\psi \colon \norm(F) \xrightarrow{\sim} (\Omega^1_{\bP^1}(D)) \simeq \cO_{\bP^1}(n-2)$ for an invertible sheaf $F$ on $C_s$ and this isomorphism can be extended to torsion-free rank one sheaves on $C_s$.
The data $\psi$ gives us the natural morphism $\cY \ra \overline{\cJ}^{n-2}_C$, which is an isomorphism. 
\end{proof}

\subsection{Calculation of the cohomology over $\cM_H$}\label{proof_of_prop_1}
In this section, we assume that $n = 5$. In this case, $V \simeq \bC^2, S = \bC^2 \setminus \{0\}$.
Moreover, the general spectral curve $C_s$ is a smooth curve of genus two branched over six points $\{ t_1, \dots, t_5, \tau_s\}$ on $\bP^1$.
A singular spectral curve occurs when the sixth point $\tau_s$ coincides with one of the five other points.
This leads to a nodal curve $C_s$ of genus two, whose desingularization $\widetilde{C}_s$ is an elliptic curve branched over $\{t_1, \dots, t_5 \} \setminus \{ \tau_s \}$ and $C_s$ can be obtained identifying two points $\tilde{\tau}^+_s$ and $\tilde{\tau}^-_s$ of $\widetilde{C}_s$.
We denote $\tilde{\tau}_s \in C_s$ the preimage of $\tau_s \in \bP^1$ by  $\pi_s \colon C_s \ra \bP^1$.
We also denote $\tilde{t}_i \in C_s$ the preimage of $t_i \in \bP^1$ by $\pi_s$ for each $i = 1, \dots, 5$.

Let us identify $\cY$ with $\overline{\cJ}^3_C$ and fix $s \in S$.
Denote by $\overline{\sigma}_s \colon \overline{\Pic}^3(C_s) \ra \overline{\Pic}^3(C_s)$ the involution defined by $F \mapsto \sigma_s^*F$.
Recall that for $x \in \bP^1$, we denote by $\xi_x$ the bundle on $\overline{\cM}$ whose fiber at $(L, \nabla, \varphi; \e \in E)$ is $L_x$.

For $y \in C_s$, let $\zeta_y$ be the torsion-free sheaf of generic rank one on $\overline{\Pic}^3(C_s)$ whose fiber over $F$ is $F_y$, and $\zeta_{\sigma(y)}$ be the the torsion-free sheaf of generic rank one on $\overline{\Pic}^3(C_s)$ whose fiber over $F$ is $F_{\sigma_s(y)}$.

Suppose $x \in \bP^1 \setminus \{t_1, \dots, t_5,\tau_s \}$ and $\pi_s^{-1}(x) = \{y^+, y^- \}$.
Moreover, suppose that $(L, \theta, \varphi) \in \cY$ corresponds to $(s, F) \in \overline{\cJ}^3_C$.
Then $(\pi^*(\xi_x))_{(L, \theta, \varphi)} = L_x = F_{y^+} \oplus F_{y^-} = (\zeta_{y^+})_{(L, \theta, \varphi)} \oplus (\zeta_{y^-})_{(L, \theta, \varphi)} $.

For $x = t_i$ and $\pi_s^{-1}(x) = \tilde{t}_i$,
we have a natural injection $(\zeta_{\tilde{t}_i})_{(L, \theta, \varphi)} \ra (\pi^*(\xi_x))_{(L, \theta, \varphi)}$.
Its cokernel is isomorphic to $(\zeta_{\tilde{t}_i})_{(L, \theta, \varphi)}$.

Similary, for $x = \tau_s$, $\pi_s^{-1}(x) =  \tilde{\tau}_s$, 
we also have a natural injection $(\zeta_{\tilde{\tau}_s})_{(L, \theta, \varphi)} \ra (\pi^*(\xi_x))_{(L, \theta, \varphi)}$ and its cokernel is isomorphic to $(\zeta_{\tilde{\tau}_s})_{(L, \theta, \varphi)}$.

\begin{Lem}
 \begin{itemize}
   \item[(i)] $\zeta_y^* = \zeta_{\sigma(y)} = \overline{\sigma}_s^*\zeta_y$ for $y \in C_s$; 
   \item[(ii)] $\zeta_{y} \not\simeq \zeta_{y'}$ for $y \neq y'$, $y, y' \in C_s$.
 \end{itemize}
\end{Lem}
\begin{proof}
(i) Since $\bigwedge^2 \pi^*(\xi_x) = \cO_{\cY}$ and $\zeta_{y^-} = \zeta_{\sigma(y^+)}$, we have the statement.

(ii)For $s \in S$, fix $y_0 \in C_s$ which is not a nodal singular point.
Consider the torsion free sheaf on $C_s \times \overline{\Pic}^3(C_s)$
whose fiber over $(y, F)$ is $F_y \otimes (F_{y_0})^{-1}$.
This sheaf is a pull back of Poincar\'{e} sheaf over $\overline{\Pic}^0(C_s) \times \overline{\Pic}^3(C_s)$ via $AJ \times \id$.
Here,
\begin{equation*}
AJ \colon C_s \ra \overline{\Pic}^0(C_s) \ :\  y \mapsto \cO_{C_s}(y - y_0) := \mathcal{H}om(\mathcal{I}_y, \cO_{C_s}(-y_0))
\end{equation*}
is the Abel-Jacobi map.
This torsion free sheaf can be viewed as a universal $C_s$ family of torsion-free sheaves of degree $0$ on $\overline{\Pic}^3(C_s)$.
In particular, two different sheaves in this $C_s$ family are not isomorphic.
Hence $F_{y} \otimes (F_{y_0})^{-1} \not\simeq F_{y'} \otimes (F_{y_0})^{-1}$ for any $y, y' \in C_s$, $y \neq y'$ and in particular $F_{y} \not\simeq F_{y'}$.
\end{proof}

Let us consider the Hitchin map, which corresponds to the projection map $p \colon \overline{\cJ}^3_C \ra S$;
\begin{equation*}
 p \colon \cY \ra S  \ : \ (L, \theta, \varphi) \mapsto \det(\theta).
\end{equation*}
This map descends to the map
$h \colon \cM_H \ra S/\mathbb{G}_m \simeq \bP^1$.
Here, $a \in \mathbb{G}_m$ acts on $s \in S$ by multiplication by $a^2$.

\begin{Prop}\label{computation_along_boundary}
Suppose $x_1, \dots, x_4 \in \bP^1$ and $x_i \neq x_j$ for $i \neq j$. Then 
\begin{equation*}
 H^i(\cM_H, \xi_{x_1}\otimes \xi_{x_2} \otimes \xi_{x_3} \otimes \xi_{x_4} \otimes (\cE | _{\cM_H})^{\otimes k}) = 0,
\end{equation*}
for any $i, k$.
\end{Prop} 
\begin{proof}
We will show this vanishing along each fiber $h^{-1}(b)$ over $b = [s] \in \bP^1$.
That is, 
\begin{equation*}
 H^i(h^{-1}(b), \xi_{x_1}\otimes \xi_{x_2} \otimes \xi_{x_3} \otimes \xi_{x_4} \otimes (\cE | _{\cM_H})^{\otimes k}|_{h^{-1}(b)}) = 0.
\end{equation*}
It is enough to prove that 
\begin{equation*}
 H^i(\overline{\Pic}^3(C_s) , \pi^*(\xi_{x_1}\otimes \xi_{x_2} \otimes \xi_{x_3} \otimes \xi_{x_4} \otimes (\cE | _{\cM_H})^{\otimes k}|_{h^{-1}(b)})) = 0.
\end{equation*}
Since $\pi^*(\cE | _{\cM_H}) = \cO_{\cY}$, we must prove 
\begin{equation}\label{restricted_cohomology}
 H^i(\overline{\Pic}^3(C_s) , \pi^*(\xi_{x_1}) \otimes \pi^*(\xi_{x_2}) \otimes \pi^*(\xi_{x_3}) \otimes \pi^*(\xi_{x_4})|_{\overline{\Pic}^3(C_s)}) = 0.
\end{equation}
But $\pi^*(\xi_{x_1}) \otimes \pi^*(\xi_{x_2}) \otimes \pi^*(\xi_{x_3}) \otimes \pi^*(\xi_{x_4})|_{\overline{\Pic}^3(C_s)}$ has a filtration with quotients $\zeta_{y_1^{\pm}} \otimes \zeta_{y_2^{\pm}} \otimes \zeta_{y_3^{\pm}} \otimes \zeta_{y_4^{\pm}}$, where 
$\pi_s(y^{\pm}_i) = x_i$.

Corollary \zcref{jacobian_vanishing_of_torsion_free_sheaves} verifies that $H^i(\overline{\Pic}^3(C_s) ,\zeta_{y_1^{\pm}} \otimes \zeta_{y_2^{\pm}} \otimes \zeta_{y_3^{\pm}} \otimes \zeta_{y_4^{\pm}}) = 0$ as a cohomology of twisted compactified Jacobians (see also Corollary \zcref{cohomology_degree0Jac}).
Therefore, the proof of this proposition follows from the proper base change theorem.
\end{proof}

\section{Bundles $\xi_x$ on $\overline{\mathcal{M}}$}\label{M_bar}
In this section, we compute the cohomology of vector bundles over $\overline{\cM}$, and prove Proposition \zcref{pre_compactified_vanish} as Proposition \zcref{computation_compactified}.  
\subsection{Parabolic $\epsilon$-connections}
We will define several additional structures on $\epsilon$-connections.
Let $l_i \subset L_{t_i}$ be a one-dimensional linear subspace for each $i = 1, \dots, n$, and $\bl = \{l_1, \dots, l_n\}$.

\begin{Def}\label{parabolic_epsilon}
A {\rm\bf {parabolic $\epsilon$-connection}} is a tuple $(L, \nabla, \varphi, \bl; \epsilon \in E)$ such that
$(L, \nabla, \varphi; \epsilon \in E)$ satisfies (1), (2), (4) of Definition \zcref{epsilon}, and the following condition
\begin{itemize}
 \item[(3')] The map $\res_{t_i}(\nabla) \colon L_{t_i} \ra L_{t_i} \otimes E$ satisfies\\
                 $\res_{t_i}(\nabla)|_{l_i} = \epsilon \nu_i$ and $(\res_{t_i}(\nabla) + \epsilon \nu_i)(L_{t_i}) \subset l_i$.
\end{itemize}
\end{Def}
We call $\bl = \{l_1, \dots, l_n \}$ {\it{parabolic structure}} and each $l_i$ {\it{parabolic direction}}.

Let $\widetilde{\overline{\cM}}$ be the moduli stack of parabolic $\epsilon$-connections. 
Denote by $\widetilde{\cM}_H \subset \widetilde{\overline{\cM}}$ the closed substack defined by the equation $\epsilon = 0$.

As explained in \cite[Remark 1]{A01}, Definition \zcref{epsilon} and \zcref{parabolic_epsilon} are equivalent in $n = 4$ case.
However, this is not the case with $n \geq 5$ as follows:
Firstly, if $\epsilon = 1$, we can set $l_i := \ker (\res_{t_i}(\nabla) - \nu_i) \subset L_{t_i}$, and so two definitions are equivalent.
Therefore, we have $\widetilde{\overline{\cM}} \setminus \widetilde{\cM}_H \simeq \overline{\cM} \setminus \cM_H = \cM$.

Next, suppose $\epsilon = 0$. If $(L, \nabla, \varphi)$ is nowhere-holomorphic, i.e., $\res_t(\nabla) \neq 0$ for every $t \in \{t_1, \dots, t_n \}$, then the parabolic structure is determined by the kernel of the residual part.

Now assume that $(L, \nabla, \varphi)$ is holomorphic at $t \in \{t_1, \dots, t_n \}$, i.e., $\res_t(\nabla) = 0$.
In this case, we cannot determine the parabolic structure uniquely.
In $n = 5$ case, the forgetful map $\widetilde{\cM}_H \ra \cM_H$ is the blowing-up along the locus formed by $(L, \nabla, \varphi)$ which is holomorphic at some point $t \in \{t_1, \dots, t_5 \}$ (see \cite[Lemma 5.2]{FL23b}).

Suppose that $(L, \nabla, \varphi, \bl; \e \in E)$ is a parabolic $\e$-connection.
Set $L' := \{ s \in L | s(t_n) \in l_n \} \subset L$, where $l_n \in \bl$.
It corresponds to the lower modification $\Elm^-_{t_n}$ as an $\e$-connection (for $\Elm^-_{t_n}$, see, e.g,  \cite[Section 2.2]{LS}).
In our $n = 5$ case, 
\begin{itemize}
 \item $L' \simeq \cO_{\bP^1} \oplus \cO_{\bP^1}(-1)$ and $\dim H^0(\bP^1, L') = 1$ if $L \simeq \cO_{\bP^1} \oplus \cO_{\bP^1}$ or $L \simeq \cO_{\bP^1}(1) \oplus \cO_{\bP^1}(-1)$ with $l_5 \not\subset \cO_{\bP^1}(1)_{t_5}$, 
 \item $L' \simeq \cO_{\bP^1}(1) \oplus \cO_{\bP^1}(-2)$ and $\dim H^0(\bP^1, L') = 2$ if $L \simeq \cO_{\bP^1}(1) \oplus \cO_{\bP^1}(-1)$ with $l_5 \subset \cO_{\bP^1}(1)_{t_5}$.
\end{itemize}

\begin{Def}
A {\rm\bf{parabolic $\e$-connection with a twisted cyclic vector}} is a tuple $(L, \nabla, \varphi, \bl, [\sigma]; \e \in E)$ such that
\begin{itemize}
 \item[(1)]$(L, \nabla, \varphi, \bl; \e \in E)$ is a parabolic $\e$-connection,
 \item[(2)]$[\sigma] \subset H^0(\bP^1, L')$ is a one-dimensional subspace generated by a nonzero section $\sigma \in H^0(\bP^1, L')$.
\end{itemize}
\end{Def}

We call a nonzero section $\sigma \in H^0(\bP^1, L')$ {\it{twisted cyclic vector}}.
Let $\widehat{\overline{\cM}}$ be the moduli stack of parabolic $\e$-connections with a twisted cyclic vector.

\begin{Lem}\label{jumping_lemma}
Suppose $n = 5$. Then, the forgetful map
 \begin{equation*}
  f \colon \widehat{\overline{\cM}} \ra \widetilde{\overline{\cM}},
 \end{equation*}
 which forgets a one-dimensional subspace generated by a twisted cyclic vector, is the blowing up along the locus $\textbf{J}$ formed by $(L, \nabla, \varphi, \bl; \e \in E)$ with bundle type $L\simeq \cO_{\bP^1}(1) \oplus \cO_{\bP^1}(-1)$ with $l_5 \subset \cO_{\bP^1}(1)_{t_5}$.
\end{Lem}
\begin{proof}
If $(L, \nabla, \varphi, \bl; \e \in E) \not\in \textbf{J}$, then $\dim H^0(\bP^1, L') = 1$, and therefore $[\sigma]$ is uniquely determined as $[\sigma] = H^0(\bP^1, L')$ with any nonzero section $\sigma \in H^0(\bP^1, L')$.

If $(L, \nabla, \varphi, \bl; \e \in E) \in \textbf{J}$, then $\dim H^0(\bP^1, L') = 2$.
In this case, each point of $\bP H^0(\bP^1, L') \simeq \bP^1$ gives us a one-dimensional subspace $[\sigma]$ of $H^0(\bP^1, L')$ with a nonzero section $\sigma \in H^0(\bP^1, L')$.
This finishes the proof of the lemma.
\end{proof}
Let us call $\textbf{J}$ {\it{jumping locus}}.
Denote by $\widehat{\cM}_H \subset \widehat{\overline{\cM}}$ the closed substack defined by the equation $\e = 0$.
We also have the open substack $\widehat{\cM} := \widehat{\overline{\cM}} \setminus \widehat{\cM}_H$, which parametrizes $\bnu$-$\mathfrak{s}l_2$-parabolic connections with a twisted cyclic vector $[\sigma] \subset H^0(\bP^1, L')$.
In the case $n = 5$, $\widehat{\cM}$ and $\widehat{\cM}_H$ are also blowing-ups of $\cM$ and $\widetilde{\cM}_H$ along their jumping loci.

\subsection{Degree $-1$ case}
In this section, we will define the moduli stacks corresponding to the case of degree $-1$. 
To calculate the cohomology of a certain vector bundle, we need a detailed description of the geometric structure of $\overline{\mathcal{M}}$ and its coarse moduli space $\overline{M}$. We will investigate these structures by focusing on the degree $-1$ case in the following sections.

Set $\bnu' =\{ \nu_i^{\pm}\}$, where $\nu_i^{\pm} = \pm \nu_i$ for $i = 1, \dots, n-1$ and $\nu_n^+ = \nu_n, \nu_n^- = 1-\nu_n$.
\begin{Def}
A {\rm\bf{$\bnu'$-$\mathfrak{s}l_2$-parabolic connection of degree $-1$}} is a triple $(L', \nabla', \varphi')$ such that
 \item[(1)] $L'$ is a rank $2$ vector bundle on $\bP^1$ of degree $-1$,
 \item[(2)] $\nabla' \colon L' \ra L' \otimes \Omega_{\bP^1}(D)$ is a connection,
 \item[(3)] $\varphi' \colon \bigwedge^2 L' \xrightarrow{\sim}  \cO_{\bP^1}(-1)$ is a horizontal isomorphism,
 \item[(4)] the residue $\res_{t_i}(\nabla')$ of the connection $\nabla'$ at $t_i$ has eigenvalues $\{\nu^+_i, \nu^-_i \}$ for each $i$ $(1 \leq i \leq n)$.
\end{Def}

As before, $E$ is a one-dimensional vector space, and $\epsilon \in E$, $L'$ is a rank $2$ vector bundle on $\bP^1$ of degree $-1$, $\nabla' \colon L' \ra L' \otimes \Omega_{\bP^1}(D) \otimes E$ is a $\bC$-linear map, and $\varphi' \colon \bigwedge^2 L' \xrightarrow{\sim} \cO_{\bP^1}(-1)$.

\begin{Def}
A {\rm\bf{parabolic $\epsilon$-connection of degree $-1$}} is a tuple $(L', \nabla', \varphi', \bl'; \epsilon \in E)$ such that
\begin{itemize}
 \item[(1)] $\nabla' (fs) = f \nabla' s + s \otimes df \otimes \epsilon$ for $f \in \cO_{\bP^1}, s \in L'$,
 \item[(2)] $\varphi' (\nabla s_1 \wedge s_2) + \varphi'(s_1 \wedge \nabla s_2) = d(\varphi'(s_1 \wedge s_2))$ for $s_1, s_2 \in L'$,
 \item[(3)] The map $\res_{t_i}(\nabla') \colon L'_{t_i} \ra L'_{t_i} \otimes E$ satisfies\\
                  $\res_{t_i}(\nabla')|_{l'_i} = \epsilon \nu_i^+$ and $(\res_{t_i}(\nabla') - \epsilon \nu_i^-)(L'_{t_i}) \subset l'_i$ for $i = 1, \dots, n$, 
 \item[(4)] $(L', \nabla')$ is irreducible.
\end{itemize}
\end{Def} 

Let $\widetilde{\overline{\cM'}}$ be the moduli stack of parabolic $\epsilon$-connections of degree $-1$.
Denote by $\widetilde{\cM}_H'$ the closed substack defined by the equation $\epsilon = 0$.
As before, the open substack $\cM' := \widetilde{\overline{\cM'}} \setminus \widetilde{\cM}_H'$ parametrizes all $\bnu'$-$\mathfrak{s}l_2$-parabolic connection of degree $-1$.

\begin{Def}\label{with_cyc_vec_odd}
A {\rm\bf{parabolic $\e$-connection of degree $-1$ with a cyclic vector}} is a tuple $(L', \nabla', \varphi', \bl', [\sigma]; \e \in E)$ such that
\begin{itemize}
 \item[(1)]$(L', \nabla', \varphi', \bl'; \e \in E)$ is a parabolic $\e$-connection of degree $-1$,
 \item[(2)]$[\sigma] \subset H^0(\bP^1, L')$ is a one-dimensional subspace generated by a nonzero section $\sigma \in H^0(\bP^1, L')$.
\end{itemize}
\end{Def}
Note that in this case, we do not need `twist' because we are already considering degree $-1$ bundles.
Let $\widehat{\overline{\cM'}}$ be the moduli stack of parabolic $\epsilon$-connections of degree $-1$ with a cyclic vector, $\widehat{\cM}_H'$ be the closed substack defined by the equation $\epsilon = 0$, and $\widehat{\cM'} := \widehat{\cM}_H' \setminus \widehat{\overline{\cM'}}$.
Then, $\widehat{\cM'}$ is the moduli stack of $\bnu'$-$\mathfrak{s}l_2$-parabolic connections of degree $-1$ with a cyclic vector.
Denote by $\widehat{M'}$ the corresponding coarse moduli space.

In the case $n = 5$, $\widehat{\overline{\cM'}}$ (respectively, $\widehat{\cM'}$, $\widehat{\cM}_H'$) is the blowing-up of $\widetilde{\overline{\cM'}}$ (respectively, $\cM'$, $\widetilde{\cM}_H'$) along the jumping locus formed by $(L', \nabla', \varphi', \bl'; \e \in E)$ with bundle type $L' \simeq \cO_{\bP^1}(1) \oplus \cO_{\bP^1}(-2)$ by using the same argument of Lemma \zcref{jumping_lemma}.

\subsection{Smooth compactification of $\widehat{M'}$}\label{lagrangian_fibrations}
We summarize here the results of Loray and Saito \cite{LS}.
In their paper, they studied the two Lagrangian fibrations on $\widehat{\cM'}$, and especially, constructed a smooth compactification of it in the case $n = 5$.

\subsubsection{Moduli space of generic connections and the two Lagrangian fibrations}\label{LSreview}
A {\textit{quasi-parabolic ${\mathfrak{s}l}_2$-bundle}} $(L, \varphi, {\bm l})$ on $(\mathbb{P}^1, D)$, ${\bm l} = \{l_1, \dots, l_n \}$, consists of a rank two vector bundle $L$ on $\mathbb{P}^1$, $\varphi: \bigwedge^2 L \xrightarrow{\sim} \cO_{\bP^1}$, and for each $i = 1, \dots, n$, a one-dimensional linear subspace $l_i \subset L_{t_i}$.

Let us introduce a notion of stability for quasi-parabolic ${\mathfrak{s}l}_2$-bundles.
For this, fix weights ${\bm w} = (w_1, \dots, w_n) \in [0, 1]^n$.
Then for any rank one subbundle $L_0 \subset L$, define the ${\bm w}$-stability index of $L_0$ to be the real number
\begin{equation}
\Stab_{{\bm w}}(L_0) := \deg(L) - 2 \deg(L_0) + \sum_{l_i \neq L_0|_{t_i}} w_i - \sum_{l_i = L_0|_{t_i}} w_i.
\end{equation}
A quasi-parabolic ${\mathfrak{s}l}_2$-bundle $(L, \varphi, {\bm l})$ is {\textit{${\bm w}$-stable}} (respectively, {\textit{${\bm w}$-semistable}}) if for every rank one subbundle $L_0 \subset L$, we have
$\Stab_{{\bm w}}(L_0) > 0$ (respectively, $\Stab_{{\bm w}}(L_0) \geq 0$).
A {\textit{parabolic ${\mathfrak{s}l}_2$-bundle}} is a quasi-parabolic ${\mathfrak{s}l}_2$-bundle together with a weight ${\bm w}$.
We say that a parabolic ${\mathfrak{s}l}_2$-bundle is {\textit{${\bm w}$-(semi)stable}} if the corresponding quasi-parabolic ${\mathfrak{s}l}_2$-bundle is ${\bm w}$-(semi)stable.

A quasi-parabolic ${\mathfrak{s}l}_2$-bundle $(L, \varphi, {\bm l})$ is {\textit{indecomposable}} if there does not exist decomposition $L = L_1 \oplus L_2$ such that each parabolic direction $l_i$ is contained either in $L_1$ or $L_2$.
It is known that a quasi-parabolic ${\mathfrak{s}l}_2$-bundle $(L, \varphi, {\bm l})$ is indecomposable if and only if it is ${\bm w}$-stable for a convenient choice of weights ${\bm w}$ (\cite[Proposition 3.4]{LS}).
Let us denote $P$ the coarse moduli space of indecomposable quasi-parabolic ${\mathfrak{s}l}_2$-bundles, and $P^{\bm w}$ the coarse moduli space of ${\bm w}$-semistable parabolic ${\mathfrak{s}l}_2$-bundles.
Then, from \cite[Proposition 3.6]{LS}, we have
\begin{equation*}
 P = \bigcup_{i, \text{finite}}P^{{\bm w}_i},
\end{equation*}
that is, $P$ can be covered by finite number of $P^{{\bm w}_i}$.
As observed in \cite{AL97, LS}, the coarse moduli space $P$ is a nonseparated scheme.

From now on, we consider the degree $-1$ case because the two Lagrangian fibrations naturally occur.
That is, a quasi-parabolic ${\mathfrak{s}l}_2$-bundle $(L', \varphi', {\bm l'})$ of degree $-1$, ${\bm l'} = \{l'_1, \dots, l'_n \}$ consists of a rank two vector bundle $L'$ on $\mathbb{P}^1$, $\varphi': \bigwedge^2 L' \xrightarrow{\sim} \cO_{\bP^1}(-1)$, and for each $i = 1, \dots, n$, a one-dimensional linear subspace $l'_i \subset L'_{t_i}$.
We define stability and indecomposability in the same way as before.

We fix here the {\textit{democratic weights}}
\begin{equation}\label{democratic}
{\bm w_0} = (w, \dots, w) \in [0, 1]^n\ \mbox{with}\ \frac{1}{n} < w < \frac{1}{n-2}.
\end{equation}
By \cite[Proposition 3.7]{LS}, for the weights ${\bm w_0} = (w, \dots, w)$ in  \zcref{democratic}, the coarse moduli space $P^{{\bm w_0}}_{-1}$ of ${\bm w_0}$-semistable parabolic ${\mathfrak{s}l}_2$-bundles of degree $-1$ is isomorphic to $\mathbb{P}H^0(\bP^1, \cO_{\bP^1}(-1) \otimes \Omega_{\bP^1}^1(D))^{\vee} \simeq \bP^{n-3}$, and consists of $(L', \varphi', {\bm l'})$ satisfying the conditions;
\begin{itemize}
 \item $L' = \cO_{\bP^1} \oplus \cO_{\bP^1}(-1)$,
 \item $l'_i \not\subset \cO_{\bP^1}$ for $i = 1, \dots, n,$
 \item not all $l'_i$ lie in the same $\cO_{\bP^1}(-1) \hookrightarrow L'$.
\end{itemize}

We also consider the coarse moduli space $M'$ of $\bnu'$-$\mathfrak{s}l_2$-parabolic connections $(L', \nabla', \varphi')$ of degree $-1$.
Note that the lower modification $\Elm_{t_n}^-$ gives us an isomorphism $M \xrightarrow{\sim} M'$, where $M$ is the coarse moduli space of $\cM$.

\begin{Rem}
\rm{By considering the ${\bm w}$-stability condition only for $\nabla$-invariant rank one subbundles, we can also define ${\bm w}$-stable parabolic $\epsilon$-connections (and so ${\bm w}$-stable $\bnu$-$\mathfrak{s}l_2$-parabolic connections).
Note that since our parabolic $\epsilon$-connections are all irreducible, every such connections is ${{\bm w}}$-stable.
However,  a parabolic $\epsilon$-connection $(L, \nabla, \varphi, \bl; \epsilon \in E)$ may have the underlying parabolic ${\mathfrak{s}l}_2$-bundle $(L, \varphi, {\bm l})$ that is not ${\bm w}$-stable.}
\end{Rem}

Now we introduce the following open subset of the moduli space $M'$.

\begin{Def}
For the democratic weight ${\bm w_0}$ in \zcref{democratic}, let us define the open subset
\begin{equation*}
 M^{{\bm w_0}}_{-1} = \{(L', \nabla', \varphi') \in M'\ | \ (L', \varphi', {\bm l'}) \in P^{{\bm w_0}}_{-1}   \}
\end{equation*}
of $M'$, which we call the moduli space of {\rm\bf{generic ${\bnu'}$-$\mathfrak{s}l_2$-parabolic connections}}.
Here, we set $l'_i := \Ker(\res_{t_i}(\nabla') - \nu_i^+) \subset L'_{t_i}$.
\end{Def}

We can define two natural Lagrangian maps on $M^{{\bm w_0}}_{-1}$.
The first one 
\begin{equation}\label{apparent_map}
\App : M^{{\bm w_0}}_{-1} \ra \bP H^0(\bP^1, \cO_{\bP^1}(-1) \otimes \Omega_{\bP^1}^1(D)) \simeq |\cO_{\bP^1}(n-3)| \simeq \bP^{n-3}_{{\bf a}}
\end{equation}
is obtained by taking the {\textit{apparent singular points}} with respect to the cyclic vector $\sigma \in H^0(\bP^1, L')$.
Here, $\bP^{n-3}_{\bf a}$ has the homogeneous coordinates ${\bf a} = (a_0 : \cdots : a_{n-3})$.
More precisely, each connection $\nabla$ on $L' = \cO_{\bP^1} \oplus \cO_{\bP^1}(-1)$ defines a $\cO_{\bP^1}$-linear map
\begin{equation*}
\cO_{\bP^1} \xrightarrow{\sigma} L' \xrightarrow{\nabla'} L' \otimes \Omega_{\bP^1}^1(D) \ra (L'/\cO_{\bP^1}) \otimes \Omega_{\bP^1}^1(D) \simeq \cO_{\bP^1}(-1) \otimes \Omega_{\bP^1}^1(D),
\end{equation*}
where the last arrow is the quotient by the subbundle defined by $\cO_{\bP^1} \hookrightarrow L'$, that is, $\nabla$ defines a map
\begin{equation*}
\varphi_{\nabla'} : \cO_{\bP^1} \ra \cO_{\bP^1}(-1) \otimes \Omega_{\bP^1}^1(D).
\end{equation*}
Its zero divisor, $\Div(\varphi_{\nabla'}) = q_1 + \cdots + q_{n-3}$, is an element of the linear system $\bP H^0(\bP^1, \cO_{\bP^1}(-1) \otimes \Omega_{\bP^1}^1(D)) \simeq |\cO_{\bP^1}(n-3)|$.
This map extends as a rational map
\begin{equation*}
\App : M' \dashrightarrow |\cO_{\bP^1}(n-3)|
\end{equation*}
on the whole moduli space with indeterminacy points on the jumping locus.

The second Lagrangian map
\begin{equation*}
\Bun : M^{{\bm w_0}}_{-1} \ra P^{{\bm w_0}}_{-1} \simeq \bP H^0(\bP^1, \cO_{\bP^1}(-1) \otimes \Omega_{\bP^1}^1(D))^{\vee} \simeq (\bP_{\bf a}^{n-3})^{\vee} \simeq \bP^{n-3}_{\bf b}
\end{equation*}
comes from the forgetful map toward the coarse moduli space $P_{-1}$ of indecomposable quasi-parabolic ${\mathfrak{s}l}_2$-bundles of degree $-1$,
\begin{equation*}
\Bun : M^{{\bm w_0}}_{-1} \ra P_{-1} ; \ (L', \nabla', \varphi') \mapsto (L', \varphi, {\bm l'})
\end{equation*}
that we restrict to the open projective chart $P^{{\bm w_0}}_{-1}$ explained above.
Here, $\bP^{n-3}_{\bf b}$ has the dual coordinates ${\bf b} = (b_0 : \cdots : b_{n-3})$.

There is a nice compactification $\overline{M^{{\bm w_0}}_{-1}}$ of the moduli space $M^{{\bm w_0}}_{-1}$ of generic ${\bnu'}$-$\mathfrak{s}l_2$-parabolic connections given by parabolic $\epsilon$-connections (for more detail, see \cite[Section 4.2]{LS}).
Set $M^{\bm w_0}_H := \overline{M^{{\bm w_0}}_{-1}} \setminus M^{{\bm w_0}}_{-1}$.
This $M^{\bm w_0}_H$ is the coarse moduli space of parabolic Higgs bundles $(L', \nabla', \varphi',  {\bm l'}; 0 \in \mathbb{C})$ such that $(L', \varphi', {\bm l'}) \in P^{{\bm w_0}}_{-1}$ and $\nabla'$ may satisfy $\det \nabla' = 0$.

\begin{Thm}[{\cite[Theorem 4.3]{LS}}]\label{LS_Thm_4_2}
Under the assumption that $\sum_i \nu_i^- \neq 0$, the morphism
\begin{equation*}
 \App \times \Bun : \overline{M^{{\bm w_0}}_{-1}} \xrightarrow{\sim} |\cO_{\bP^1}(n -3)| \times |\cO_{\bP^1}(n -3)|^{\vee} \simeq \bP_{\bf a}^{n-3} \times \bP_{\bf b}^{n-3}
\end{equation*}
is an isomorphism.
Moreover, by restriction, we also obtain the isomorphism
\begin{equation*}
 \App \times \Bun|_{M^{\bm w_0}_H} : M^{\bm w_0}_H \xrightarrow{\sim} \Sigma, 
\end{equation*}
where $\Sigma$ is the incidence variety for the duality. \qed
\end{Thm}

\subsubsection{Blowing-ups of $\bP_{\bf a}^{2} \times \bP_{\bf b}^{2}$ as a smooth compactification of $\widehat{M'}$}\label{blowing_ups}
Until the end of this chapter, assume that $n = 5$.
As explained in \cite[Section 6.2]{LS}, in $n = 5$ case, we can construct a smooth compactification of the full coarse moduli space $\widehat{M'}$ corresponding to $\widehat{\cM'}$ by blowing-up $\overline{M^{{\bm w_0}}_{-1}} \simeq \bP_{\bf a}^{2} \times \bP_{\bf b}^{2}$ appropriately.
This leads that $\overline{M}$ and $M_H$, the coarse moduli spaces corresponding to $\overline{\cM}$ and $\cM_H$, and their blowing-ups are all rational projective schemes.

Firstly, let us specify some important sets in $\bP_{\bf a}^{2} \times \bP_{\bf b}^{2}$ by using the coordinates ${\bf a}$ and ${\bf b}$.
In $\bP^2_{\bf a}$, let us denote by $\Delta$ the image of the diagonal though the map
\begin{equation*}
 \Sym \colon \bP^1 \times \bP^1 \ra \Sym^2 \bP^1 = \bP^2_{\bf a}\ : \ (q_1, q_2) \mapsto (z-q_1)(z-q_2),
\end{equation*}
that is the conic $\Delta \colon \{ a_1^2 -4a_0a_2 = 0\}$, which corresponds to the locus of double roots $q_1=q_2$.
It is naturally parametrized by the base curve
\begin{equation*}
 \bP^1 \ra \Delta \ : \ q \mapsto (q^2:-2q:1).
\end{equation*}
The locus $q = t_i$ of poles give us five special points on the conic $\Delta$, that is $(a_0 : a_1: a_2) = (t_i^2: -2t_i:1)$, and we will denote by $\Delta_i : \{t_i^2 a_2 + t_i a_1+ a_0 = 0 \}$ the line tangent to $\Delta$ at this point.
Any two of those lines intersect at a point $\Delta_i \cap \Delta_j = \{\Delta_{i,j} \}$ with a coordinate $(a_0 : a_1: a_2) = (t_it_j : -t_i - t_j : 1)$.

In $\bP^2_{\bf b}$, the dual of $\bP^2_{\bf a}$, we have the dual conic $\Pi := \Delta^* : \{b_1^2 -b_0b_2 = 0 \}$, which is also parametrized by the base curve
\begin{equation*}
 \bP^1 \ra \Pi \ : \ z \mapsto (1 : z : z^2).
\end{equation*}
We get five points $D_i := \Delta^*_i$ on the dual conic $\Pi$ defined by $(b_0 : b_1 : b_2) = (1 : t_i : t_i ^2)$, and $10$ lines $\Pi_{i, j} := \Delta_{i, j}^*$ passing though both $D_i$ and $D_j$ with the equation $\{t_i t_j b_0 - (t_i + t_j)b_1 + b_2 = 0\}$.
See also \cite[Figure 3]{LS}.

Let $\Sigma \subset \bP^2_{\bf a} \times \bP^2_{\bf b}$ the incidence variety defined by $\{a_0b_0 + a_1b_1 + a_2b_2 = 0 \}$.
The conic $\Delta \subset \bP^2_{\bf a}$ lifts-up as a rational curve $\Gamma \subset \Sigma$ parametrized by 
\begin{equation*}
 \bP^1 \ra \Gamma \ : \ q \mapsto ((1: -2q : q^2), (q^2: q: 1)).
\end{equation*}
It is defined by the equations
\begin{equation*}
a_1^2 = 4a_0a_2,\  a_0b_0 = a_2b_2,\  \text{and} \ \ 2a_2b_2+a_1b_1 = 0.
\end{equation*}
Inside of $\Sigma$, we also define five lines
\begin{equation*}
\Gamma_i := \Delta_i \times \{D_i \}
\end{equation*}
and $10$ more lines
\begin{equation*}
 \Gamma_{i, j} := \{\Delta_{i, j}\} \times \Pi_{i, j}.
\end{equation*}

Next, we will specify the $16$ $(-1)$-curves in the chart $\mathcal{S}$ that corresponds to the moduli space $P_{-1}^{\bm w}$ with weights ${\bm w} = (w, \dots, w)$ with $1/3 < w < 3/5$.
The chart $\mathcal{S}$ is isomorphic to the blowing-up of the five points $D_i$ in $\bP^2_{\bf b}$, and by definition, it is the del Pezzo surface of degree $4$.
Let us denote $\Pi_i$ the exceptional divisors corresponding to this blowing-up, and we keep the notation $\Pi_{i, j}$ and $\Pi$ for the strict transformations of them.
These curves $\Pi_i, \Pi_{i, j}$ and $\Pi$ constitute $16$ $(-1)$-curves on $\mathcal{S}$.
See \zcref{16del} for the list of the corresponding parabolic ${\mathfrak{s}l}_2$-bundles of degree $-1$.

\begin{table}
 \centering
  \begin{tabular}[p]{|c|c|c|c|}
   \hline
            & $L$ &  $\{i,j,k,m,n\}=\{t_1, t_2, t_3, t_4, t_5\}$ & locus\\
   \hline \hline           
   $5$ & $\cO_{\bP^1} \oplus \cO_{\bP^1}(-1)$ & $l_i \subset \cO_{\bP^1}$ & $\Pi_i, \tilde{F}_i, G_i$\\
   \hline
    $10$ & $\cO_{\bP^1}\oplus \cO_{\bP^1}(-1)$ & $l_k, l_m, l_n \subset \cO_{\bP^1}(-1)$ & $\Pi_{i, j}, \tilde{F}_{i, j}, G_{i, j}$\\
   \hline
   $1$ & $\cO_{\bP^1}\oplus \cO_{\bP^1}(-1)$ & $l_i, l_j, l_k, l_m, l_n \subset \cO_{\bP^1}(-2)$ & $\Pi, \tilde{F}, G$\\
   \hline
  \end{tabular}
\caption{$16$ special lines on $\mathcal{S}$ in \cite[Table 1]{LS}.}
\label{16del}
\end{table}

To get a smooth compactification of the full coarse moduli space $\widehat{M'}$, we have to blow up $16$ curves $\Gamma_i, \Gamma_{i, j}$ and $\Gamma$.
More precisely, we have to blow up;
\begin{itemize}
 \item[(1)] $\Gamma_i$ for those connections on a bundle having the parabolic $l'_i \in \cO_{\bP^1}$,
 \item[(2)] $\Gamma_{i, j}$ for those connections on a bundle having $l'_i, l'_j \in \cO_{\bP^1}$,
 \item[(3)] $\Gamma$ for those connections on $\cO_{\bP^1}(1) \oplus \cO_{\bP^1}(-2)$ with $[\sigma]$.
\end{itemize}
Note that case (2) has been forgotten to be mentioned in \cite{LS}.
We will explain the blowing-up procedure in case (1).
The computations in cases (2) and (3) are similar to those of case (1).

For simplicity, set $D = 0 + 1 + t_1 + t_2 + \infty$, and $\kappa_i := \nu_i^+ - \nu_i^-$, for $i = 0, 1, t_1, t_2, \infty$.
Firstly, we consider a curve in $\bP^2_{\ba} \times \bP^2_{\bf b}$ formed by $\bnu'$-${\mathfrak{sl}}_2$-parabolic connections $\{(L'_t, \nabla'_t, \varphi'_t, \bl'_t)\}$ of degree $-1$ such that,
\begin{itemize}
 \item if $t \neq 0$, then $\App \times \Bun(L'_t, \nabla'_t, \varphi'_t, \bl'_t) \in \bP^2_{\ba} \times \bP^2_{\bf b} \setminus \Sigma \ (\simeq M^{w_0}_{-1}$),
 \item if $t = 0$, then $L'_0 \simeq \cO_{\bP^1} \oplus \cO_{\bP^1}(-1)$, $(l'_{t_1})_0 \in \cO_{\bP^1}$, and other $(l'_i)_0$ are generic (especially $(l'_{\infty})_0 \not\in \cO_{\bP^1})$.
\end{itemize}
By applying the elementary transformation $\Elm^+_{\infty}$, the corresponding family of connections $\{(L_t, \nabla_t, \varphi_t, \bl_t) \}$ is described as follows:
$L_t \simeq \cO_{\bP^1} \oplus \cO_{\bP^1}$, $\varphi_t \colon \bigwedge^2 L_t \simeq \cO_{\bP^1}$ for all $t$, 
\begin{equation*}
    (l_0)_t= \begin{pmatrix}1\\ 0\end{pmatrix}, (l_1)_t= \begin{pmatrix}1\\ 1\end{pmatrix}, (l_{t_1})_t= \begin{pmatrix}1\\ u_1^t\end{pmatrix}, (l_{t_2})_t= \begin{pmatrix}1\\ u_2^t\end{pmatrix},  (l_{\infty})_t= \begin{pmatrix}0\\ 1\end{pmatrix},
\end{equation*}
and
\begin{equation*}
 \nabla_t = \nabla_0 + c_1^t \Theta_1(u_1^t, u_2^t) + c_2^t \Theta_2(u_1^t, u_2^t),
\end{equation*}
where 
\begin{equation*}
\nabla_0 :=d+\begin{pmatrix}\nu_0^-&0\\\rho&\nu_0^+\end{pmatrix}\frac{dz}{z}
+\begin{pmatrix}\nu_1^--\rho&\kappa_1 +\rho\\ -\rho&\nu_1^++\rho\end{pmatrix}\frac{dz}{z-1}
+\sum_{i=1}^{2}\begin{pmatrix}\nu_{t_i}^-&\kappa_{t_i}u_i\\ 0&\nu_{t_i}^+\end{pmatrix}\frac{dz}{z-t_i},
\end{equation*}
with $\rho= \nu_{0}^-+ \nu_{1}^-+ \nu_{\infty}^-+ \nu_{t_1}^- + \nu_{t_2}^-$,
and
\begin{equation*}
\Theta_i:=\begin{pmatrix}0&0\\1-u_i&0\end{pmatrix}\frac{dz}{z}
+\begin{pmatrix}u_i&-u_i\\ u_i&-u_i\end{pmatrix}\frac{dz}{z-1}
+\begin{pmatrix}-u_i&u_i^2\\-1&u_i\end{pmatrix}\frac{dz}{z-t_i},
\end{equation*}
\begin{equation*}
 c_1^t=-t\kappa_{t_1}+t^2\cdot c_1,\ \ \ c_2^t=c_2,\ \ \ u_1^t=\frac{1}{t}, \ \ \ u_2^t=u_2 \ \ (c_i, u_i \in \bC).
 \end{equation*}

The explicit calculations in \cite[Section 6.3]{LS} tell us that, after applying $\Elm_{\infty}^-$, the limiting point of this curve when $t \ra 0$ tends to $\Gamma_{t_1} (\subset \Sigma)$ given by
\begin{equation*}
 (a_0 : a_1: a_2) \ra (t_1 q_2 : -t_1 - q_2 : 1), \ (b_0 : b_1 : b_2) \ra (1 : t_1 : t_1^2),
\end{equation*}
with apparent points given by
\begin{equation*}
q_1 = t_1, \ q_2 = \frac{t_2(c_2(u_2-1)- \rho -\kappa_{t_1})}{c_2(u_2-t_2)-\rho-\kappa_{t_1}}.
\end{equation*}
To distinguish the connections having the same limiting point, we have to blow up $\Gamma_{t_1}$.
 Let us denote $F_{t_1}$ the exceptional divisor.
 $F_{t_1}$ is isomorphic to a $\bP^2$-bundle over $\bP^1$ parametrized by $q_2 \in \Delta_{t_1}, \tilde{u}_2 \in \Pi_{t_1}$, and the additional parameter $c$ that corresponds to the endomorphisms of the underlying parabolic ${\mathfrak{s}l}_2$-bundles such that $(\tilde{u}_2, c) \in \bP^2$.

To get the desired connections, we have to blow-up $F_{t_1}$ again because we do not have a parameter corresponding to $c_1$.
Precisely, we have to blow-up the surface defined in $F_{t_1}$ by
\begin{equation*}
 (\rho+\kappa_{t_1})u+\kappa_{t_1}t_1(t_1+q_2)v-\kappa_{t_1}t_1q_2w=0,
\end{equation*}
with the fiber coordinates $(u: v: w) \in \bP^2$.
Let $E_{t_1}$ be the exceptional divisor, and $\tilde{F}_{t_1}$ be the strict transformation of $F_{t_1}$.
Then, we get the remaining parameter corresponding to $c_1$, and $E_{t_1}$ parametrizes the connections with the condition $l'_{t_1} \in \cO_{\bP^1}$. 

In the same way, we can recover the other connections on a bundle having the parabolic $l'_i \in \cO_{\bP^1}$ for $i = 0, 1, t_2, \infty$.

In the case (2) (respectively, (3)), we have to blow-up $\Gamma_{i, j}$ (respectively, $\Gamma$) twice.
Denote by $F_{i, j}$ (respectively, $F$) the first exceptional divisor, $E_{i, j}$ (respectively, $E$) the second exceptional divisor, and $\tilde{F}_{i, j}$ (respectively, $\tilde{F}$) the strict transformation of $F_{i, j}$ (respectively, $F$). 
For more detail, see \cite[Section 6.4]{LS}.

Let $\widehat{\bP_{\bf a}^2 \times \bP_{\bf b}^2} \ra \bP_{\bf a}^2 \times \bP_{\bf b}^2$ denote the corresponding blowing-up along the 16 curves $\{\Gamma_i, \Gamma_{i, j}, \Gamma \}$ that we blow up twice described as above.
Then, we get a birational map
\begin{equation*}
  \widetilde{\App \times \Bun} \colon \widehat{\overline{M'}} \dashrightarrow \widehat{\bP_{\bf a}^2 \times \bP_{\bf b}^2}
\end{equation*}
such that, by restriction, we obtain an injective morphism
\begin{equation*}
 \widetilde{\App \times \Bun}|_{\widehat{M'}} \colon \widehat{M'} \hookrightarrow \widehat{\bP_{\bf a}^2 \times \bP_{\bf b}^2}.
\end{equation*}
That is, we get a smooth compactification of $\widehat{M'}$.
\subsection{The indeterminacy locus}\label{rationality_section}
In this section, we will study the indeterminacy locus of the birational map
\begin{equation*}
 \widetilde{\App \times \Bun} \colon \widehat{\overline{M'}} \dashrightarrow \widehat{\bP_{\bf a}^2 \times \bP_{\bf b}^2}.
\end{equation*}

Firstly, we will define apparent maps from each $\widehat{\overline{\cM}}$ and $\widehat{\overline{\cM'}}$ to $\bP^2_{\bf a}$.
Let $(L, \nabla, \varphi, \bl, [\sigma]; \e \in E) \in \widehat{\overline{\cM}}$ be a parabolic $\e$-connection with a twisted cyclic vector $[\sigma] \subset H^0(\bP^1, L')$, $L' := \{ s \in L | s(t_5) \in l_5 \} \subset L$.
Then $\nabla \sigma \in H^0(\bP^1, L' \otimes \Omega^1_{\bP^1}(D) \otimes E)$, and so 
\begin{equation*}
\nabla \sigma \wedge \sigma \in H^0(\bP^1, \Omega^1_{\bP^1}(D - t_5) \otimes E), \nabla \sigma \wedge \sigma \neq 0.
\end{equation*}
Let $\{q_1, q_2 \}$ be zeros of $\nabla \sigma \wedge \sigma$. 
Here $q_i \in \bP^1$, and we get the next map
\begin{equation*}
\App \colon \widehat{\overline{\cM}} \ra \bP^2_{\ba} = \Sym^2(\bP^1)\ :\  (L, \nabla, \varphi, \bl, [\sigma]; \e \in E) \mapsto [q_1, q_2].
\end{equation*}

Next, let $(L', \nabla', \varphi', \bl', [\sigma]; \e \in E) \in \widehat{\overline{\cM'}}$ be a parabolic $\e$-connection of degree $-1$ with a cyclic vector $[\sigma] \subset H^0(\bP^1, L')$.
Then we can also define the map
\begin{equation*}
\App \colon \widehat{\overline{\cM'}} \ra \bP^2_{\ba} = \Sym^2(\bP^1)\ :\  (L', \nabla', \varphi', \bl', [\sigma]; \e \in E) \mapsto [q_1, q_2].
\end{equation*}
in the same way.

Note that we have a commutative diagram
\begin{equation}\label{comm_diag_elm}
\xymatrix{
\widehat{\overline{\cM}}\ar[rr]^-{\Elm_{t_5}^-}_-{\sim}\ar[dr]_-{\App}&\ar@{}[d]|{\circlearrowright}&\widehat{\overline{\cM'}}\ar[dl]^-{\App}\\
&\bP^2_{\bf a}&
}
\end{equation}

Let us consider the Hitchin map
\begin{equation*}
 \hat{h} \colon \widehat{\cM}'_H \ra \bP^1 \ ; \ (L', \nabla', \varphi', \bl', [\sigma] ; 0 \in E) \mapsto [\det(\nabla')].
\end{equation*}
For each $b = [s] \in \bP^1$, the fiber $\hat{h}^{-1}(b)$ has $16$ $2$-torsion points $\{p_i \}_{i = 1}^{16}$ corresponding to the classes $[F] \in \hat{h}^{-1}(b)$  such that $F \simeq \sigma_s^* F$ (cf. \cite[Section 5.3]{DP}).
Varying $b \in \bP^1$, we get $16$ lines $\{P_i \}_{i = 1}^{16}$ in $\widehat{\cM}'_H$.
These lines correspond to the set of parabolic $\epsilon$-connections that admit unstable underlying parabolic bundles.
Especially, we can classify all of them.

\begin{Prop}\label{classification_of_P}
Suppose $(L', \nabla', \varphi', \bl', [\sigma] ; 0 \in E) \in \cup_{i = 1}^{16} P_i$.
Then the underlying parabolic ${\mathfrak{s}l}_2$-bundle $(L', \varphi', \bl')$ is one of the following types:
\begin{itemize}
     \item[(I)] $L' = \cO_{\bP^1} \oplus \cO_{\bP^1}(-1)$, $\cO_{\bP^1}$ contains $1$ parabolic direction and $\cO_{\bP^1}(-1)$ contains $4$ parabolic directions;
     \item[(I\hspace{-1.2pt}I)] $L' = \cO_{\bP^1} \oplus \cO_{\bP^1}(-1)$, $\cO_{\bP^1}$ contains $2$ parabolic directions and $\cO_{\bP^1}(-1)$ contains $3$ parabolic directions;
     \item[(I\hspace{-1.2pt}I\hspace{-1.2pt}I)] $L' = \cO_{\bP^1}(1) \oplus \cO_{\bP^1}(-2)$, $\cO_{\bP^1}(1)$ contains no parabolic direction and $\cO_{\bP^1}(-2)$ contains all parabolic directions.
    \end{itemize}
Moreover, the apparent singular points $\{q_1, q_2\}$ of the corresponding $(\nabla', [\sigma])$ satisfy the following conditions:
 \begin{itemize}
  \item[(I)]  One of $\{q_1, q_2\}$ equals to $t_i$, $t_i \in \{t_1, \dots, t_5 \}$, i.e., $[q_1, q_2] \in \Delta_i$;
  \item[(I\hspace{-1.2pt}I)] $\{q_1, q_2 \} = \{t_i, t_j \}$, $t_i, t_j \in \{t_1, \dots, t_5 \}$, and $i \neq j$, i.e., $[q_1, q_2] \in \Delta_{i, j}$;
  \item[(I\hspace{-1.2pt}I\hspace{-1.2pt}I)] $q_1 = q_2$, i.e., $[q_1, q_2] \in \Delta$.
 \end{itemize}    
\end{Prop}
\begin{proof}
Suppose $(L', \nabla', \varphi', \bl', [\sigma] ; 0 \in E) \in \cup_{i = 1}^{16} P_i$.
Firstly, we will study the underlying parabolic ${\mathfrak{s}l}_2$-bundle $(L', \varphi', \bl')$.
Since it corresponds to a $2$-torsion point explained above, $\nabla'$ satisfies the condition $(-1)\cdot \nabla' \sim \nabla'$, that is, there exists $P \in GL_2(\bC)$ such that $(-1) \cdot \nabla' = P^{-1} \nabla'  P$.
Then, the underlying parabolic bundle $(L', \bl')$ admits a non-trivial automorphism, and therefore, $(L', \varphi', \bl')$ is unstable with respect to any weights ${\bm w} \in [0, 1]^5$.
Especially, $(L', \varphi', \bl')$ is unstable with respect to the weight ${\bm w}_c = (1/2, \dots, 1/2)$.

Let $Y^{{\bm w}_c}_{-1}$ be the coarse moduli space of ${\bm w}_c$-semistable parabolic ${\mathfrak{sl}}_2$-Higgs bundles of degree $-1$ over $(\bP^1, D)$, where ${\bm w}_c = (1/2, \dots, 1/2)$. 
Assume $(L', \theta', \varphi', \bl') \in Y^{{\bm w}_c}_{-1}$ and $(L', \varphi', \bl')$ is ${\bm w}_c$-unstable.
Then, by the same argument of \cite[Corollary 3.2]{FL23b} in the degree $-1$ case, we can conclude that there are only three possibilities (I), (I\hspace{-1.2pt}I), and (I\hspace{-1.2pt}I\hspace{-1.2pt}I) in the statement.
This implies that there are exactly $16$ ${\bm w}_c$-unstable parabolic ${\mathfrak{s}l}_2$-bundles that admits a ${\bm w}_c$-semistable Higgs field $\theta'$, see \zcref{16}.
The group $(\mathbb{Z}/2 \mathbb{Z})^4$ acts transitively on it as elementary transformations (\cite[section 2.2]{FL23b}).
Therefore, our $(L', \varphi', \bl')$ should be one of them, because there are only $16$ $2$-torsion points on each fiber of the Hitchin map $\hat{h} \colon \widehat{\cM}'_H \ra \bP^1$.

Next, we will check the apparent singular points of $(\nabla', [\sigma])$.
For simplicity, suppose $D = 0 + 1 + \infty + t_1 + t_2$, $L' = \cO_{\bP^1} \oplus \cO_{\bP^1}(-1)$, and the parabolic directions $l_0$ over $0$ lie in $\cO_{\bP^1}$ and $l_1, l_{\infty}, l_{t_1}, l_{t_2}$ over $1, \infty, t_1, t_2$ lie in $\cO_{\bP^1}(-1)$.
Then, following the same argument of \cite[Corollary 3.2, Remark 3.3]{FL23b} in the degree $-1$ case, we can also check that any ${\bm w}_c$-semistable Higgs field on it is of the form
\begin{eqnarray}\label{betagamma}
\theta'=\left(
\begin{array}{ccc} 
0 & \beta  \\
\gamma & 0 \\
\end{array}
\right)
\end{eqnarray}
with 
\begin{displaymath}
\left\{ \begin{array}{ll}
\beta: \cO_{\bP^1}(-1) \ra \Omega_{\bP^1}(0), \quad \beta \neq 0\\
\gamma: \cO_{\bP^1} \ra \cO_{\bP^1}(-1) \otimes \Omega_{\bP^1}(1 + t_1 + t_2 + \infty)\;, \quad \gamma\neq 0\,.
\end{array} \right.
\end{displaymath}
Especially, $\gamma$ vanishes at $\{0\}$.
This implies that one of the apparent singular points $\{q_1, q_2 \}$ of $\theta'$ is $\{0 \}$.
Since $\mathbb{G}_m$-action on $(L', \theta', \varphi', \bl')$ does not change the apparent singular points, our $\nabla'$ also has the same apparent singular points.

Conversely, such $\nabla'$ satisfies the equation
\begin{eqnarray*}
(-1) \cdot \left(
\begin{array}{ccc} 
0 & \beta  \\
\gamma & 0 \\
\end{array}
\right)
=
\left(
\begin{array}{ccc} 
1 & 0  \\
0 & -1 \\
\end{array}
\right)^{-1}
\left(
\begin{array}{ccc} 
0 & \beta  \\
\gamma & 0 \\
\end{array}
\right)
\left(
\begin{array}{ccc} 
1 & 0  \\
0 & -1 \\
\end{array}
\right),
\end{eqnarray*}
so $(L', \nabla', \varphi', \bl' ; 0 \in E)$ corresponds to a $2$-torsion point of $\hat{h}^{-1}(b)$ with $b = [-\beta \gamma] \in \bP^1$.

Any other Higgs bundle admitting a ${\bm w}_c$-unstable parabolic ${\mathfrak{s}l}_2$-bundle can be obtained by performing an elementary transformation, and therefore, the proposition is proven.
\end{proof}

\begin{table}
 \centering
  \begin{tabular}[p]{|c|c|c|c|}
   \hline
            & $L$ &  $\{i,j,k,m, n\}=\{t_1, t_2, t_3, t_4, t_5\}$ & locus\\
   \hline \hline           
   $5$ & $\cO_{\bP^1}\oplus \cO_{\bP^1}(-1)$ & $l_i\subset \cO_{\bP^1}$ and $l_j, l_k, l_m, l_n\subset \cO_{\bP^1}(-1)$ & $P_i$\\
   \hline
   $10$ & $\cO_{\bP^1} \oplus \cO_{\bP^1}(-1)$ & $l_i, l_j \subset \cO_{\bP^1}$ and $l_k, l_m l_n \subset \cO_{\bP^1}(-1)$ & $P_{i, j}$\\
   \hline
   $1$ & $\cO_{\bP^1}(1)\oplus \cO_{\bP^1}(-2)$ & $l_i, l_j, l_k, l_m, l_n \subset \cO_{\bP^1}(-2)$ & $P$ \\
   \hline
  \end{tabular}
\caption{$16$ unstable parabolic bundles admitting semistable Higgs fields.}
\label{16}
\end{table}

We will introduce the convenient notation following \cite{DP}.
Let $\Odd$ be the set of all subsets of $\{t_1, t_2, t_3, t_4, t_5 \}$ of odd cardinality.
Then $\{\Delta_i, \Delta_{i, j}, \Delta \}$, $\{\Pi_i, \Pi_{i, j}, \Pi \}$, $\{\Gamma_i, \Gamma_{i, j}, \Gamma \}$, $\{F_i, F_{i, j}, F \}$, $\{\tilde{F}_i, \tilde{F}_{i, j}, \tilde{F} \}$, and the 16 lines $\{P_i \}_{i = 1}^{16}$ can be naturally labeled by the subset $I \in \Odd$ as follows:
\begin{itemize}
 \item if $I = \{ t_i\}$, then we will set $\{\Delta_I, \Pi_I, \Gamma_I, F_I, \tilde{F}_I \}$ to be $\{\Delta_i, \Pi_i, \Gamma_i, F_i, \tilde{F}_i \}$, and $P_I \subset \widehat{\cM}'_H$ to be the type (I) of the Proposition \zcref{classification_of_P} with $l_i \in \cO_{\bP^1}$;
 \item if $\sharp I = 3$, then we will set $\{\Delta_I, \Pi_I, \Gamma_I, F_I, \tilde{F}_I \}$ to be $\{\Delta_{i, j}, \Pi_{i, j}, \Gamma_{i, j}, F_{i, j}, \tilde{F}_{i, j} \}$, and $P_I \subset \widehat{\cM}'_H$ to be the type (I\hspace{-1.2pt}I) of the Proposition \zcref{classification_of_P} with $l_i, l_j \in \cO_{\bP^1}$ for  $i, j \not\in I$;
 \item if $I = \{t_1, t_2, t_3, t_4, t_5 \}$, then we will set $\{\Delta_I, \Pi_I, \Gamma_I, F_I, \tilde{F}_I \}$ to be $\{\Delta, \Pi, \Gamma, F, \tilde{F} \}$, and $P_I \subset \widehat{\cM}'_H$ to be the type (I\hspace{-1.2pt}I\hspace{-1.2pt}I) of the Proposition \zcref{classification_of_P}.
\end{itemize}

From the discussion in \zcref{blowing_ups}, we get the following birational map
\begin{equation*}
 \widetilde{\App \times \Bun} \colon \widehat{\overline{M'}} \dashrightarrow \widehat{\bP_{\bf a}^2 \times \bP_{\bf b}^2}
\end{equation*}
such that $\widetilde{\App \times \Bun}(\widehat{\overline{M'}} \setminus \cup_I P_I) = \widehat{\bP_{\bf a}^2 \times \bP_{\bf b}^2} \setminus \cup_I \tilde{F}_I$.
Let $\widehat{\overline{M''}}$ denote the blowing-up of $\widehat{\overline{M'}}$ along the above 16 lines $\{P_I\}_{I \in \Odd}$, and $\{G_I \}_{I \in \Odd}$ denote the exceptional divisors.
This blowing-up corresponds to considering the opposite order extensions of the underlying parabolic ${\mathfrak{s}l}_2$-bundle $(L', \varphi', \bl')$.
Therefore, $(L', \varphi', \bl')$ becomes ${\bm w}_c$-stable, and by comparing \zcref{16del} and \zcref{16}, we can check that $(L', \varphi', \bl') \in \Pi_I$.
Here, we still have the endomorphism $\phi$ of $(L', \bl')$ such that 
\begin{eqnarray}\label{phi}
\phi \sim
\left(
\begin{array}{ccc} 
c & 0  \\
0 & 0 \\
\end{array}
\right)
\end{eqnarray}
with $c \in \bC$.
Therefore, we get the following lemma:
\begin{Lem}\label{point_of_G_I}
The point of $G_I$ corresponds to an isomorphic class of a tuple $(L', \nabla', \varphi', \bl', [\sigma], \phi; 0 \in E)$ such that 
\begin{itemize}
 \item $(L', \nabla', \varphi', \bl', [\sigma]; 0 \in E)$ is a parabolic $\epsilon$-connection of degree $-1$ with a cyclic vector with $\epsilon = 0$;
 \item the underlying parabolic bundle $(L', \varphi', \bl') \in \Pi_I$;
 \item the apparent singular points $[q_1, q_2] \in \Delta_I$;
 \item $\phi \in \End(L', \bl')$ that satisfies the condition \zcref{phi}.
\end{itemize}
\end{Lem}\qed

We will show that the birational map $\widetilde{\App \times \Bun}$ extends to an isomorphism. 

\begin{Thm}\label{extended_app_bun}
The birational map $ \widetilde{\App \times \Bun} \colon \widehat{\overline{M'}} \dashrightarrow \widehat{\bP_{\bf a}^2 \times \bP_{\bf b}^2}$ extends to an isomorphism
 \begin{equation*}
  \widehat{\App \times \Bun} \colon \widehat{\overline{M''}} \xrightarrow{\sim} \widehat{\bP_{\bf a}^2 \times \bP_{\bf b}^2}.
 \end{equation*}
\end{Thm}
\begin{proof}
We will show that the locus $\tilde{F}_I$ corresponds to the exceptional divisor $G_I$ for every $I \in \Odd$.
Since there is a $(\mathbb{Z}/2\mathbb{Z})^4$-symmetry as elementary transformations, we only need to check the correspondence between $\tilde{F}_{t_1}$ and $G_{t_1}$ under the assumption $D = 0 + 1 + t_1 + t_2 + \infty$.

As explained in \zcref{blowing_ups}, the locus $\tilde{F}_{t_1}$ is isomorphic to a $\bP^2$-bundle over $\bP^1$ parametrized by $q_2 \in \Delta_{t_1}$, $\tilde{u}_2 \in \Pi_{t_1}$ and the additional parameter $c$ such that $(\tilde{u}_2, c) \in \bP^2$. 
On the other hand, by the Lemma \zcref{point_of_G_I}, we can construct the unique element $(L', \nabla', \varphi', \bl', [\sigma], \phi; 0 \in E) \in G_{t_1}$ from each point $(q_2, \tilde{u}_2, c) \in \tilde{F}_{t_1}$.
Therefore, we can extend the map $\widetilde{\App \times \Bun}$ to the morphism $\widehat{\App \times \Bun}$ such that $\widehat{\App \times \Bun}(G_{t_1}) \simeq \tilde{F}_{t_1}$, concluding the proof of the theorem.
\end{proof}

\subsection{Calculation of the cohomology over $\overline{\cM}$}\label{proof_of_prop_2}
Set $\widehat{\cE} := \cO_{\widehat{\overline{\cM}}}(\widehat{\cM}_H)$, $\widehat{\cE'} := \cO_{\widehat{\overline{\cM'}}}(\widehat{\cM}'_H)$, and $\widehat{\cE''} := \cO_{\widehat{\overline{\cM''}}}(\widehat{\cM}''_H)$, where $\widehat{\cM}''_H := \widehat{\overline{\cM''}} \setminus \widehat{\cM'}$.

\begin{Lem}\label{higher_direct_image_of_E} 
\begin{itemize}
\item[(1)] For any integers $q \geq 0$ and $k = 1, -2, -3$, we have
\begin{equation}
 R^q \App _* (\widehat{\cE}^{\otimes k}) \simeq H^q(\bP^2, \cO_{\bP^2}(k)) \otimes \cO_{\bP^2}(k).
\end{equation}
\item[(2)]
\begin{equation}\label{direct_image_structure}
 R^q \App_* \cO_{\widehat{\overline{\cM}}} \simeq \begin{cases}
                                                                                  \cO_{\bP^2}, & q = 0,\\
                                                                                  0, & q \neq 0.
\end{cases}                                                                                  
\end{equation}
\end{itemize}
\end{Lem}
\begin{proof}
(1) Since we have the commutative diagram \zcref{comm_diag_elm}, it is enough to consider the degree $-1$ case. 
 From the Theorem \zcref{extended_app_bun}, there is a diagram
\begin{equation}\label{comm_diag_App}
  \begin{CD}
       \widehat{\overline{\cM''}} @>f_1>> \widehat{\overline{\cM'}}\\
      @Vf_2VV    @VV\App V   \\
     \bP^2_{\bf a} \times \bP^2_{\bf b} @>pr_1>> \bP^2_{\bf a}
  \end{CD}
\end{equation}
where $f_1$ is a blowing-up of $\widehat{\overline{\cM'}}$ along the $16$ lines $\{ P_I \}_{I \in \Odd}$ and $f_2$ is a composition of morphisms
\begin{equation*}
 \widehat{\overline{\cM''}} \ra \widehat{\overline{M''}} \ra \widehat{\bP^2_{\bf a} \times \bP^2_{\bf b}} \ra \bP^2_{\bf a} \times \bP^2_{\bf b}.
\end{equation*}
On the other hand, since $k = 1, -2, -3$, we have
\begin{equation*}
 R^q f_{1, *} (\widehat{\cE''})^{\otimes k} \simeq \begin{cases}
                                                                                  (\widehat{\cE'})^{\otimes k}, & q = 0,\\
                                                                                  0, & q \neq 0,
                                                                                 \end{cases}
\end{equation*}
\begin{equation*}
 R^q f_{2, *} (\widehat{\cE''})^{\otimes k} \simeq \begin{cases}
                                                                                  \cO_{\bP^2_{\bf a} \times \bP^2_{\bf b}}(k \Sigma), & q = 0,\\
                                                                                  0, & q \neq 0.
                                                                                 \end{cases}
\end{equation*}
From the commutativity of the above diagram \zcref{comm_diag_App} and the Grothendieck spectral sequence,
\begin{equation*}
 \begin{split}
  R^q\App_* (\widehat{\cE'})^{\otimes k} &\simeq R^q \App_* f_{1, *} (\widehat{\cE''})^{\otimes k}\\
                                                                &\simeq R^q(\App \circ f_1)_* (\widehat{\cE''})^{\otimes k}\\
                                                                &\simeq R^q(pr_1 \circ f_2)_*(\widehat{\cE''})^{\otimes k}\\
                                                                &\simeq R^q pr_{1, *} f_{2, *}(\widehat{\cE''})^{\otimes k}\\
                                                                &\simeq R^q pr_{1, *} \cO_{\bP^2_{\bf a} \times \bP^2_{\bf b}}(k \Sigma)\\
                                                                &\simeq H^q(\bP_{\bf b}^2, \cO_{\bP_{\bf b}^2}(k)) \otimes \cO_{\bP_{\bf a}^2}(k).
 \end{split}
\end{equation*}
We can prove (2) in the same way. 
\end{proof}

\begin{Prop}\label{cohomology_of_compactification}
We have
\begin{itemize}
\item[(1)]
  $H^i(\widehat{\overline{\mathcal{M}}}, \mathcal{O}_{\widehat{\overline{\mathcal{M}}}}) =
   \begin{cases}
    \mathbb{C}, & i = 0\\
             0, & i \neq 0,
   \end{cases}$
\item[(2)]
  $H^i(\widehat{\cM}_H, \mathcal{O}_{\widehat{\cM}_H}) =
   \begin{cases}
    \mathbb{C}, & i = 0\\
             0, & i \neq 0.
   \end{cases}$
\end{itemize}
\end{Prop}
\begin{proof}
(1) Denote by $\overline{M}$ the coarse moduli space corresponding to $\overline{\mathcal{M}}$.
From the discussion in \zcref{rationality_section}, we know that $\overline{M}$ is a rational projective scheme. 
Since $\widehat{\overline{M}}$ is the blowing-up of $\overline{M}$, 
\begin{equation*}
 H^i(\widehat{\overline{\mathcal{M}}}, \mathcal{O}_{\widehat{\overline{\mathcal{M}}}}) = H^i(\widehat{\overline{M}}, \mathcal{O}_{\widehat{\overline{M}}}) =
  \begin{cases}
   \mathbb{C}, & i = 0\\
            0, & i \neq 0.
  \end{cases}
 \end{equation*} 
We can prove (2) in the same way.
\end{proof}

Denote by $\xi_+$ the line bundle on $\widehat{\overline{\cM}}$ whose fiber over $(L, \nabla, \varphi, \bl, [\sigma]; \e \in E)$ equals $[\sigma] \subset H^0(\bP^1, L')$.
\begin{Lem}\label{filtration_vanish}
For $i \geq 0$,
 \begin{itemize}
  \item[(1)]$H^i(\widehat{\overline{\cM}}, \cO_{\widehat{\overline{\cM}}}(-\widehat{\cM}_H)) = 0$,
  \item[(2)]$H^i(\widehat{\overline{\cM}}, (\xi_+)^{\otimes 2}(-\widehat{\cM}_H)) = 0$,
  \item[(3)]$H^i(\widehat{\overline{\cM}}, ((\xi_+)^*)^{\otimes 2}(-\widehat{\cM}_H)) = 0$,
  \item[(4)]$H^i(\widehat{\overline{\cM}}, (\xi_+)^{\otimes 4}(-\widehat{\cM}_H)) = 0$,
  \item[(5)]$H^i(\widehat{\overline{\cM}}, ((\xi_+)^*)^{\otimes 4}(-\widehat{\cM}_H)) = 0$.
 \end{itemize}
\end{Lem}
\begin{proof}
(1) Consider the exact sequence
\begin{equation*}
0 \ra \cO_{\widehat{\overline{\cM}}}(-\widehat{\cM}_H) \ra \cO_{\widehat{\overline{\cM}}} \ra \cO_{\widehat{\overline{\cM}}}/\cO_{\widehat{\overline{\cM}}}(-\widehat{\cM}_H) \ra 0.
\end{equation*}
From Proposition \zcref{cohomology_of_compactification}, the natural map
\begin{equation*}
H^i(\widehat{\overline{\cM}}, \cO_{\widehat{\overline{\cM}}}) \ra H^i(\widehat{\overline{\cM}}, \cO_{\widehat{\overline{\cM}}}/\cO_{\widehat{\overline{\cM}}}(-\widehat{\cM}_H)) = H^i(\widehat{\cM}_H, \cO_{\widehat{\cM}_H})
\end{equation*}
is bijective. So we have the first statement.

(2)
Fix $x_0 \in \bP^1 \setminus \{t_1, \dots, t_5 \}$, $\omega \in \Omega^1_{\bP^1}(D)_{x_0}$, $\omega \neq 0$.
The correspondence
\begin{equation*}
\sigma \mapsto (\nabla \sigma \wedge \sigma)(x_0)\omega^{-1} \in \bigwedge^2 L_{x_0} \otimes E
\end{equation*}
defines a map $(\xi_+)^{\otimes 2} \ra \widehat{\cE} = \cO_{\widehat{\overline{\cM}}}(\widehat{\cM}_H)$.
This map vanishes along the locus $\App^{-1}(l_{x_0})$ where $l_{x_0} \subset \bP^2_{\ba}$ is the line formed by $[q_1, q_2]$ such that one of $q_i$ is $x_0$.
Therefore, we get an isomorphism $(\xi_+)^{\otimes 2} \xrightarrow{\sim} \App^*(\cO_{\bP^2}(-1)) \otimes \widehat{\cE}$.
So $(\xi_+)^{\otimes 2}(- \widehat{\cM}_H) \simeq \App^*(\cO_{\bP^2}(-1))$.

Therefore, by the projection formula and \zcref{direct_image_structure}, we have
\begin{equation*}
 \begin{split}
 H^i(\widehat{\overline{\cM}}, (\xi_+)^{\otimes 2}(- \widehat{\cM}_H)) &= H^i(\widehat{\overline{\cM}}, \App^*(\cO_{\bP^2}(-1)))\\
                                                                        &= H^i(\bP^2, \cO_{\bP^2}(-1))\\
                                                                        &= 0
 \end{split}
\end{equation*}
for all $i$.

(3) From the above discussion, $(\xi_+^*)^{\otimes 2}(- \widehat{\cM}_H) \simeq \widehat{\cE}^{\otimes (-2)} \otimes \App^*(\cO_{\bP^2}(1))$. 
Therefore, by Leray spectral sequence, the projection formula, and Lemma \zcref{higher_direct_image_of_E}, we have
 \begin{equation*}
  \begin{split}
   H^i(\widehat{\overline{\cM}}, (\xi_+^*)^{\otimes 2}(- \widehat{\cM}_H)) &= H^i(\widehat{\overline{\cM}}, \widehat{\cE}^{\otimes (-2)} \otimes \App^*(\cO_{\bP^2}(1)))\\
     &= \bigoplus_{p + q = i} H^p(\bP^2, R^q\App_*(\widehat{\cE}^{\otimes (-2)} \otimes \App^*(\cO_{\bP^2}(1))))\\
     &= \bigoplus_{p + q = i} H^p(\bP^2, R^q\App_*\widehat{\cE}^{\otimes (-2)} \otimes \cO_{\bP^2}(1))\\
     &= 0.
  \end{split}
 \end{equation*}

(4) We have $(\xi_+)^{\otimes 4}(- \widehat{\cM}_H) \simeq \widehat{\cE} \otimes \App^*(\cO_{\bP^2}(-2))$.
Therefore,
 \begin{equation*}
  \begin{split}
   H^i(\widehat{\overline{\cM}}, (\xi_+)^{\otimes 4}(- \widehat{\cM}_H)) &= H^i(\widehat{\overline{\cM}}, \widehat{\cE} \otimes \App^*(\cO_{\bP^2}(-2)))\\
   &= \bigoplus_{p + q = i} H^p(\bP^2, R^q\App_*(\widehat{\cE} \otimes \App^*(\cO_{\bP^2}(-2))))\\
   &= \bigoplus_{p + q = i} H^p(\bP^2, R^q\App_*\widehat{\cE} \otimes \cO_{\bP^2}(-2))\\
   &= H^i(\bP^2, \cO_{\bP^2}(-1))^{\oplus 3}\\
   &= 0.
  \end{split}
 \end{equation*}

(5) We have $(\xi_+^*)^{\otimes 4}(- \widehat{\cM}_H) \simeq \widehat{\cE}^{\otimes (-3)} \otimes \App^*(\cO_{\bP^2}(2))$.
Therefore,
 \begin{equation*}
  \begin{split}
   H^i(\widehat{\overline{\cM}}, (\xi_+^*)^{\otimes 4}(- \widehat{\cM}_H)) &= H^i(\widehat{\overline{\cM}}, \widehat{\cE}^{\otimes (-3)} \otimes \App^*(\cO_{\bP^2}(2)))\\
   &= \bigoplus_{p + q = i} H^p(\bP^2, R^q\App_*(\widehat{\cE}^{\otimes (-3)} \otimes \App^*(\cO_{\bP^2}(2))))\\
   &= \bigoplus_{p + q = i} H^p(\bP^2, R^q\App_*\widehat{\cE}^{\otimes (-3)} \otimes \cO_{\bP^2}(2))\\
   &= H^{i-2}(\bP^2, \cO_{\bP^2}(-1))\\
   &= 0.
  \end{split}
 \end{equation*}
\end{proof}

\begin{Prop}\label{computation_compactified}
Suppose $x_1, \dots, x_4 \in \bP^1$. Then
\begin{equation*}
 H^i(\overline{\cM}, \xi_{x_1} \otimes \xi_{x_2} \otimes \xi_{x_3} \otimes \xi_{x_4}(-\cM_H)) = 0,
\end{equation*} 
for any $i$.
\end{Prop}
\begin{proof}
Let us show that
\begin{equation}\label{blowing_uped_vanishing}
H^i(\widehat{\overline{\cM}}, \xi_{x_1} \otimes \xi_{x_2} \otimes \xi_{x_3} \otimes \xi_{x_4}(-\widehat{\cM}_H)) = 0.
\end{equation}
Without loss of generality, we may assume that $x_1, \dots, x_4$ are not equal to $t_5$.
Since we have an inclusion $L' \hookrightarrow L$, the natural map $\xi_+ \ra \xi_{x_i}$ is injective and its cokernel is isomorphic to $(\xi_+)^*$.
We use this map to identify $\xi_+$ with a subbundle of $\xi_{x_i}$.
Then $\xi_{x_1} \otimes \xi_{x_2} \otimes \xi_{x_3} \otimes \xi_{x_4}$ has a filtration $\{ \cF_k \}$ with quotients $\cF_k/\cF_{k-1} = (\xi_+)^{\otimes 2}, (\xi_+)^{\otimes 4}, (\xi_+^*)^{\otimes 2}, (\xi_+^*)^{\otimes 4}$, or $\cO_{\widehat{\overline{\cM}}}$.

It follows that $H^i(\widehat{\overline{\cM}}, (\cF_k/\cF_{k-1})(-\widehat{\cM}_H)) = 0$ by Lemma \zcref{filtration_vanish}.
Therefore, we get \zcref{blowing_uped_vanishing}.
Since the forgetful map $\widehat{\overline{\cM}} \ra \overline{\cM}$ is the composition of blowing-ups, this implies the statement.
\end{proof}

\section{Cohomology of the structure sheaf of $\cM$}\label{Proof_of_Theorem2}
Suppose $n = 5$.
In this section, we will compute the cohomology of the structure sheaf of $\cM$, the moduli space of parabolic connections over $(\bP^1, t_1+\cdots + t_5)$, following the strategy given by D. Arinkin  \cite{AF} based on his discussion with R. Fedorov.
\subsection{Cohomology of compactified Jacobians}

Fix $g \geq 0$.
Let $p_C \colon C \rightarrow S$ be a family of projective integral curves with planar singularities of arithmetic genus $2$ over a base scheme $S$.
Let $\overline{J}^3_C$ be the moduli space of pairs $(s, F)$, where $s \in S$ and $F$ is a torsion-free sheaf of degree $3$ of generic rank one on $C_s$. 
Let $\overline{\mathcal{J}}^3_C$ be the $\mu_2$-gerbe over $\overline{J}^3_C$ explained in \zcref{m2_gerbe} (\cite[section 4.2]{MSY}).

We have already shown in Corollary \zcref{higher_direct_image_of_strucuture_sheaf} that
\begin{equation}\label{exterior_product}
 R^{\bullet}p^3_* (\cO_{\overline{\mathcal{J}}^3_C}) = \bigwedge^{\bullet} R^1p_{C, *} \cO_C,
\end{equation}
where $p^3\colon \overline{\mathcal{J}}^3_C \ra S$, $p_C\colon C \ra S$.
In our case, 
\begin{equation*}
S = H^0(\bP^1, \Omega_{\bP^1}^{\otimes 2}(t_1 + \cdots + t_5)) \setminus \{0 \} \simeq \bC^{2} \setminus \{0\}.
\end{equation*}
Set $H^1 := H^1(\bP^1, \cO_{\bP^1} \oplus (\Omega_{\bP^1}(t_1 + \cdots + t_5))^{-1}) = H^1(\bP^1, \pi_{s, *}\cO_{C_s}) \simeq H^1(C_s, \cO_{C_s})$, where $\pi_s : C_s \ra \bP^1$ and $C_s$ is a smooth curve.

\begin{Lem}
In the derived category $D^b(S)$, we have an isomorphism
 \begin{equation}\label{decomp}
  Rp_{C, *}\cO_{C} \simeq \cO_S \oplus H^1[-1] \otimes \cO_S. 
 \end{equation}
\end{Lem}
\begin{proof}
By using Koll\'{a}r's decomposition theorem, we have
\begin{equation*}
 Rp_{C, *}\omega_C \simeq p_{C, *}\omega_C \oplus R^1p_{C, *}\omega_C [-1].
\end{equation*}
Applying Grothendieck duality to both sides, we get
\begin{equation*}
Rp_{C, *}\cO_C \simeq p_{C, *}\cO_C \oplus R^1p_{C, *}\cO_C [-1].
\end{equation*}
Therefore, we have to show that $R^ip_{C, *}\cO_C \simeq H^i(C_s, \cO_{C_s}) \otimes \cO_S$, that is, $R^ip_{C, *}\cO_C$ are trivial bundles over $S \simeq \bC^{2} \setminus \{0\}$ for $i = 0, 1$.
Firstly, it is well-known in this case that $p_{C, *}\cO_C \simeq \cO_S$. Next, by using \zcref{exterior_product} and Matsushita's theorem (\cite{Matsushita}, \cite[Example 3.5]{MSY}), we have
\begin{equation*}
 R^1p_{C, *}\cO_C \simeq R^1p^3_* (\cO_{\overline{\mathcal{J}}^3_C}) \simeq \Omega_S \simeq H^1(C_s, \cO_{C_s}) \otimes\cO_S.
\end{equation*}
\end{proof}

\begin{Prop}\label{claim2}
We have 
 \begin{equation}
  Rp^3_*(\cO_{\overline{\mathcal{J}}^3_C}) = \Sym^{\bullet}(H^1[-1]) \otimes \cO_S.
 \end{equation}
 Here, $\Sym^{\bullet}(H^1[-1]) = \bigoplus(\bigwedge^i H^1)[-i]$.
\end{Prop}
\begin{proof}
\zcref{exterior_product} and \zcref{decomp} give us an isomorphism as graded algebras.
 We have to show that the direct image $Rp^3_*(\cO_{\overline{\mathcal{J}}^3_C})$ is formal.
 
 Recall that $Rp^3_*(\cO_{\overline{\mathcal{J}}^3_C}) \simeq R{\mathcal{H}om}_{S}(\mathcal{N}^{-3}, \mathcal{N}^{-3})$.
Here, the right-hand side depends only on the formal neighborhood of the zero section of $\cJ_C$. 
Since $\cJ_{C}$ and $\mathcal{N}^{-3}$ make sense over the entire $S$ (more precisely, over $\bC^{2}$), this means that there exists an object $\cF \in D^b(S)$ such that $H^i(\cF) \simeq (\bigwedge^i H^1) \otimes \cO_S$ and $Rp^3_*(\cO_{\overline{\mathcal{J}}^3_C}) \simeq \cF$.
  Since $S \simeq \bC^{2} \setminus \{0 \}$ and $\cF$ has locally free cohomology, $\cF$ is formal.
Therefore, $Rp^3_*(\cO_{\overline{\mathcal{J}}^3_C})$ is also formal.
\end{proof}

\subsection{Cohomology of the Hitchin Systems (Theorem \zcref{main2})}\label{section_cohomology_of_hitchin_systems}
Recall that we can identify $\overline{\cJ}^{3}_C$ with $\cY$, the moduli space of $\mathfrak{s}l_2$-Higgs bundles over $(\bP^1, D)$ by Proposition \zcref{BNR}, and we have the next commutative diagram:
\begin{equation}
  \begin{CD}
      \cY \simeq \overline{\mathcal{J}}^{3}_C  @>p^{3}>> S \simeq \bC^{2} \setminus \{ 0\}\\
      @V/\mathbb{G}_mV\pi V    @V/\mathbb{G}_mVqV   \\
   \cM_H @>h>> \bP^{1}_{(2)}.
  \end{CD}
\end{equation}
Note that $a \in \mathbb{G}_m$ acts on $s \in S$ by multiplication by $a^2$.
Therefore we get a weighted projective space $\bP^{1}_{(2)}$.
Consider the $\mathbb{G}_m$-equivariant isomorphism of Proposition \zcref{claim2} in our case:
\begin{equation}\label{equivariant_case}
 Rp^{3}_*(\cO_{\cY}) = \Sym^{\bullet}(H^1[-1]) \otimes \cO_{S}.
\end{equation}

Set $V := H^0(\bP^1, \Omega_{\bP^1}^{\otimes 2} (t_1 + \cdots + t_5)) \simeq \bC^{2}$, which is isomorphic to $(H^1)^*$ by the Serre duality.
Note that $S = V \setminus \{0\}$.
Let us take the cohomology of both sides of \zcref{equivariant_case}.
Here we consider the cohomology $H^{\bullet}(\cY, \cO_{\cY})$ as a bigraded algebra with respect to the cohomological grading and the grading by weight of $\mathbb{G}_m$.

\begin{Prop}\label{cohomology_of_Y}
There is an isomorphism of bigraded algebras
 \begin{equation*}
   H^{\bullet}(\cY, \cO_{\cY}) \simeq (\Sym^{\bullet}(V^*) \oplus \Sym^{\bullet}(V) \otimes \det V [-1]) \otimes \Sym^{\bullet}(H^1[-1]).
 \end{equation*}
\end{Prop}
\begin{proof}
It is well known that
\begin{equation*}
H^i(\mathbb{A}^2 \setminus \{0\}, \cO) = \begin{cases}
                                     \bC[x_1, x_2] & \text{if}\ i = 0,\\
                                     x_1^{-1} x_2^{-1} \cdot \bC[x_1^{-1}, x_2^{-1}] & \text{if}\ i = 1,\\
                                     0 & \text{otherwise}.
                                    \end{cases}
\end{equation*}
Therefore, the statement follows from this and \zcref{equivariant_case}.
\end{proof}

Note that the weight $-j$ component $H^{\bullet}(\cY, \cO_{\cY})_{(-j)}$ of $H^{\bullet}(\cY, \cO_{\cY})$ is equal to $H^{\bullet}(\cM_H, (\cE|_{\cM_H})^{\otimes j})$.
Therefore, $H^{\bullet}(\cY, \cO_{\cY}) = \bigoplus_{i, j} H^i(\cM_H, (\cE|_{\cM_H})^j)$.

\begin{Lem}\label{canonical_bundle_over_M_H}
 We have
  $\omega_{\cM_H} \simeq (\cE|_{\cM_H})^{\otimes (-2)}$.
\end{Lem}
\begin{proof}
By using Leray spectral sequence, we have
\begin{equation*}
 H^i(\cM_H, \cO_{\cM_H}) = \bigoplus_{p+q=i} H^q(\bP^{1}_{(2)}, R^ph_* \cO_{\cM_H}).
\end{equation*} 
From Proposition \zcref{cohomology_of_compactification}, we know
$H^q(\bP^{1}_{(2)}, R^2h_* \cO_{\cM_H}) = 0$ for all $q$.
On the other hand, since the fiber $h^{-1}(b)$ over $b = [s] \in \bP^1_{(2)}$ is isomorphic to $\overline{\Pic}^3(C_s)$, 
\begin{equation*}
 (R^2h_*\cO_{\cM_H})_b = H^2(h^{-1}(b), \cO_{h^{-1}(b)}) = \bC
\end{equation*}
for every $b \in \bP^1_{(2)}$.
Therefore $R^2h_*\cO_{\cM_H} \simeq \cO_{\bP^1_{(2)}}(-2)$, and by Grothendieck duality, $h_*\omega_{\cM_H} \simeq \cO_{\bP^1_{(2)}}(-2)$.
Since $(\cE|_{\cM_H})^{\otimes 2} \simeq h^*\cO_{\bP^1_{(2)}}(2)$, we get the statement.
\end{proof}

Denote by $\cM_{H(k)}$ the $k$-th infinitesimal neighborhood of $\cM_H$.
This means that $\cM_{H(k)} \subset \overline{\cM}$ is the closed substack defined by the sheaf of ideal $\cO_{\overline{\cM}}(-k \cM_H) \subset \cO_{\overline{\cM}}$.
From the filtration
\begin{equation*}
\cdots \subset \cO_{\overline{\cM}}(- \cM_H) \subset \cO_{\overline{\cM}} \subset \cO_{\overline{\cM}}(\cM_H) \subset \cdots,
\end{equation*}
we have a short exact sequence
\begin{equation*}
 0 \ra (\cE|_{\cM_H})^{\otimes (j - 1)} \ra (\cE|_{\cM_{H(2)}})^{\otimes j} \ra (\cE|_{\cM_H})^{\otimes j} \ra 0,
\end{equation*}
and the associated long exact sequence of cohomology
\begin{equation}\label{MH_long_exact_sq}
 \begin{split}
\cdots &\rightarrow H^i(\cM_H, (\cE|_{\cM})^{\otimes (j - 1)}) \rightarrow H^i(\cM_{H(2)}, (\cE|_{\cM_{H(2)}})^{\otimes j}) \rightarrow H^i(\cM_H, (\cE|_{\cM_H})^{\otimes j})\\
  &\xrightarrow{\delta^{i, j}} H^{i + 1}(\cM_H, (\cE|_{\cM_H})^{\otimes (j -1)}) \rightarrow H^{i+1}(\cM_{H(2)}, (\cE|_{\cM_{H(2)}})^{\otimes j}) \rightarrow  \cdots.
   \end{split}
\end{equation}

Especially, we get
\begin{equation*}
 \delta^{i, j}\colon H^i(\cM_H, (\cE|_{\cM_H})^{\otimes j}) \ra H^{i + 1}(\cM_H, (\cE|_{\cM_H})^{\otimes (j-1)}),
\end{equation*}
a differential of bidegree $(1, -1)$ on $\bigoplus_{i, j} H^i(\cM_H, (\cE|_{\cM_H})^{\otimes j})$.

From now on, $\cE|_{\cM_H} = \cO_{\overline{\cM}}(\cM_H)|_{\cM_H}$ is abbreviated to $\cE$ even over $\cM_H$ for simplicity.
Consider $\delta^{0, 2}\colon H^0(\cM_H, \cE^{\otimes 2}) \ra H^1(\cM_H, \cE)$, which corresponds to $V^* \ra H^1$.

Denote by $H^i_j$ the weight $-j$ part of $H^i(\cY, \cO_{\cY})$, i.e., $H^i_j = H^i(\cM_H, \cE^{\otimes j})$, and $H_{>0} := \oplus_{j > 0} H^{\bullet}_j$.
\begin{Prop}\label{cohomology_of_bigraded_algebra}
If $\delta^{0, 2}$ is an isomorphism, then any non-zero cohomology of the complex $(H^{\bullet}(\cY, \cO_{\cY}), \delta)$ are
$\bC$ in bidegree $(0, 0)$, and $\det(V) \otimes \det(H^1) = \bC$ in bidegree $(3, -2)$.
Especially, $(H_{> 0}, \delta)$ is exact.
\end{Prop}
\begin{proof}
First, we will show that $\delta^{0, 2}$ determines all $\delta^{\bullet, \bullet}$.
Since 
\begin{equation*}
H^2_0 = H^2(\cM_H, \cO_{\cM_H}) = 0
\end{equation*}
from Proposition \zcref{cohomology_of_compactification} (2), we have $\delta^{1, 1} = 0$.
Thus, $\delta^{0, 2}$ determines the differential $\delta^{\bullet, \bullet}$ on the subalgebra of $H^{\bullet}(\cY, \cO_{\cY})$ generated by $H^0_2, H^1_1$.
Moreover, this subalgebra corresponds to $\Sym^{\bullet}(V^*) \otimes \Sym^{\bullet}(H^1[-1])$ of Proposition \zcref{cohomology_of_Y}.

Let us show that $\delta^{0, 2}$ also determines the remaining component.
Since $\dim \cM_H = 3$, we have $H^{4}(\cM_H, \cE^{\otimes (-3)}) = 0$.
Therefore $\delta^{3, -2} = 0$.
Moreover, the product
\begin{equation*}
 H^i_j \otimes H^{3 - i}_{-2-j} \rightarrow H^{3}_{-2}
\end{equation*}
is a non-degenerate pairing because of Serre duality.
Therefore, $\delta^{0, 2}$ also determines the subalgebra generated by $H^{3}_{-4}, H^{2}_{-3}$ and it corresponds to the remaining component.

It remains to prove that $(H_{>0}, \delta)$ is exact.
From the first part of this proof, we know that 
\begin{equation*}
 H_{\geq 0} := \oplus_{j \geq 0}H_j^{\bullet} = \Sym^{\bullet}(V^*) \otimes \wedge^{\bullet} H^1.
\end{equation*}
Choose a basis $x_1, x_2$ of $V^*$, and set $y_j = \delta^{0,2}(x_j)$ for $j = 1,2$.
Since $\delta^{0,2}$ is an isomorphism, $y_1, y_2$ form a basis of $H^1$.
Thus, 
\begin{equation*}
  H_{\geq 0} = \bC[x_1, x_2] \otimes \wedge^{\bullet}\langle y_1, y_2 \rangle.
\end{equation*}
By construction, $\delta(x_j) = y_j$ and $\delta(y_j) = 0$ because $y_j \in H^1_1$ and $H^2_0 = 0$.
The graded Leibniz rule
\begin{equation}\label{graded_Leibniz}
 \delta(xy) = \delta(x) y + (-1)^i x \delta(y)
\end{equation}
gives us $\delta = \sum_{j=1}^2 y_j \wedge \frac{\partial}{\partial x_j}$.
In particular, $\delta$ decreases the polynomial degree by one and increases the exterior degree by one.
Let $y_1^*, y_2^*$ be the dual basis of $(H^1)^*$.
For each $j$, set 
\begin{equation*}
\iota_{y_j^*} \colon \wedge^{q}H^1 \ra \wedge^{q-1}H^1
\end{equation*}
be the contraction, characterized by $\iota_{y_j^*}(y_i) = \delta_{ij}$ and \zcref{graded_Leibniz}.
Extend this operator $\bC[x_1, x_2]$-linearly, and define
\begin{equation*}
 h := \sum_{j = 1}^2 x_j \iota_{y_j^*}. 
\end{equation*}
Then, we have
\begin{equation}\label{elementary_identity}
 \delta h + h \delta = (p+q)\id
\end{equation}
on $\Sym^p(V^*)\otimes \wedge^q H^1$.
Both $\delta$ and $h$ preserve $p + q$, and we have $p+q > 0$ on $H_{>0}$.
Therefore, we may define $s := \frac{h}{p+q}$ and \zcref{elementary_identity} gives us
\begin{equation*}
 \delta s + s \delta = \id.
\end{equation*}
We conclude that $(H_{>0}, \delta)$ is exact.
\end{proof}

\begin{Cor}\label{proof_of_main2}
In the hypothesis of the Proposition \zcref{cohomology_of_bigraded_algebra}, the restriction map $H^{\bullet}(\overline{\cM}, \cO_{\overline{\cM}}) \ra H^{\bullet}(\cM, \cO_{\cM})$ is an isomorphism. Therefore,
\begin{equation*}
  H^i(\cM, \cO_{\cM}) = \begin{cases}
                                     \bC & \text{if}\ i = 0,\\
                                     0 & \text{if}\ i > 0.
                                    \end{cases}
 \end{equation*}
\end{Cor}
\begin{proof}
Denote by $j \colon \cM \hookrightarrow \overline{\cM}$ and $i \colon \cM_H \hookrightarrow \overline{\cM}$ the natural embeddings.
Put $\cQ_k := \cO_{\overline{\cM}}(k \cM_H)/ \cO_{\overline{\cM}}$, $\cQ_0 := 0$, and $\cQ := \varinjlim \cQ_k$.
Their successive quotients fit into
\begin{equation}\label{successive_quotient}
 0 \ra \cQ_{k-1} \ra \cQ_k \ra i_*\cE^{\otimes k} \ra 0
\end{equation}
for $k \geq 1$.
Since $j_* \cO_{\cM} = \varinjlim \cO_{\overline{\cM}}(k \cM_H)$, and $\cQ = j_* \cO_{\cM}/\cO_{\overline{\cM}}$, we have
\begin{equation*}
0 \ra \cO_{\overline{\cM}} \ra j_* \cO_{\cM} \ra \cQ \ra 0.
\end{equation*}
Therefore, we only need to prove $H^i(\overline{\cM}, \cQ) = 0$ for $i \geq 0$.

Let $a \in H^i(\overline{\cM}, \cQ_k)$, and let $x \in H^i_k$ be its image under 
\begin{equation*}
 \varphi^{i, k} \colon H^i(\overline{\cM}, \cQ_k) \ra H^i(\cM_H, \cE^{\otimes k})
\end{equation*}
induced by \zcref{successive_quotient} so that $x = \varphi^{i, k}(a)$.
If $k \geq 2$, the existence of this lift gives us $\partial x = \partial \varphi^{i, k}(a) = 0$ with 
$ \partial \colon H^i(\cM_H, \cE^{\otimes k}) \ra H^{i+1}(\overline{\cM}, \cQ_{k-1})$. 
Moreover, the commutativity of the diagram
\begin{equation*}
\xymatrix{
{H^i(\cM_H, \cE^{\otimes k})}\ar[rr]^-{\partial}_-{}\ar[dr]_-{\delta^{i, k}}&\ar@{}[d]|{\circlearrowright}& H^{i+1}(\overline{\cM}, \cQ_{k-1}) \ar[dl]^-{\varphi^{i+1, k-1}}\\
& H^{i+1}(\cM_H, \cE^{\otimes (k-1)}) &
}
\end{equation*}
shows $\delta^{i, k}(x) = 0$.
If $k = 1$, then $\delta^{i, 1}(x) = 0$ because $H^{i+1}_0 = 0$. 

First, suppose $i > 0$.
The exactness of $(H_{>0}, \delta)$ in Proposition \zcref{cohomology_of_bigraded_algebra} gives $y \in H^{i-1}_{k+1}$ with $\delta^{i-1, k+1}(y) = x$.
Let $c := \partial y \in H^i(\overline{\cM}, \cQ_k)$ be its image under the connecting homomorphism for 
\begin{equation*}
 0 \ra \cQ_{k} \ra \cQ_{k+1} \ra i_*\cE^{\otimes k+1} \ra 0.
\end{equation*}
The commutativity of the diagram
\begin{equation*}
\xymatrix{
{H^{i-1}(\cM_H, \cE^{\otimes k+1})}\ar[rr]^-{\partial}_-{}\ar[dr]_-{\delta^{i-1, k+1}}&\ar@{}[d]|{\circlearrowright}& H^{i}(\overline{\cM}, \cQ_{k}) \ar[dl]^-{\varphi^{i, k}}\\
& H^{i}(\cM_H, \cE^{\otimes k}) &
}
\end{equation*}
identifies the image of $c$ in $H^i_k$ with $\delta^{i-1, k+1}(y) = x$.
Therefore, $\varphi^{i, k}(a-c) = 0$, and hence $a - c$ comes from $H^i(\overline{\cM}, \cQ_{k-1})$.
The exactness of the same long exact sequence says that $c$ becomes zero in $H^i(\overline{\cM}, \cQ_{k+1})$.
Thus, in the filtered colimit, $a$ is represented by a class of one lower pole order.
Repeating finitely many times reaches $\cQ_0 = 0$.

If $i = 0$, the exactness of $(H_{>0}, \delta)$ already implies $x = 0$, since there is no group in cohomological degree $-1$.
The class $a$ therefore comes directly from $H^0(\overline{\cM}, \cQ_{k-1})$ so we can use the same induction.

Therefore, we have proved that $\varinjlim H^i(\overline{\cM}, \cQ_k) = 0$ for $i \geq 0$.
Since $\overline{\cM}$ is a complete algebraic stack, the filtered colimits statement in \cite[\href{https://stacks.math.columbia.edu/tag/0GQX}{Tag 0GQX}]{stacks-project} tells us the commutativity of taking cohomology and filtered colimits.
Therefore, we have
\begin{equation*}
 \begin{split}
  H^i(\overline{\cM}, \cQ) &= H^i(\overline{\cM}, \varinjlim \cQ_k)\\
                                        &= \varinjlim H^i(\overline{\cM}, \cQ_k)\\
                                        &= 0.
 \end{split}
\end{equation*}
for $i \geq 0$.
Hence $H^{\bullet}(\overline{\cM}, \cO_{\overline{\cM}}) \ra H^{\bullet}(\cM, \cO_{\cM})$ is an isomorphism. 
\end{proof}

\begin{Prop}\label{delta_isomorphism}
$\delta^{0, 2}$ is an isomorphism.
\end{Prop}
\begin{proof}
We will sketch the proof following \cite{AF}.
The map $\delta^{0, 2}$ appears as the connecting homomorphism.
That is, 
\begin{equation*}
 \begin{split}
0 &\rightarrow H^0(\cM_H, \cE) \rightarrow H^0(\cM_{H(2)}, \cE^{\otimes 2}) \rightarrow H^0(\cM_H, \cE^{\otimes 2})\\
  &\xrightarrow{\delta^{0, 2}} H^1(\cM_H, \cE) \rightarrow H^1(\cM_{H(2)}, \cE^{\otimes 2}) \rightarrow  \cdots.
   \end{split}
\end{equation*}
Since $\dim H^0(\cM_H, \cE^{\otimes 2}) = \dim H^0(\bP^1_{(2)}, \cO_{\bP^1_{(2)}}(2)) = \dim H^0(\bP^1, \cO_{\bP^1}(1)) = 2$, and
$\dim H^1(\cM_H, \cE) = \dim H^1(h^{-1}(b), \cO_{h^{-1}(b)}) = 2$ with some $b \in \bP^1_{(2)}$, we only need to show that $\delta^{0, 2}$ is injective.
We will show an equivalence statement: For each $i > 0$, the restriction map
\begin{equation}
 H^0(\cM_{H(2)}, \cE^{\otimes i}) \ra H^0(\cM_H, \cE^{\otimes i})
\end{equation}
is zero.

Consider the second infinitesimal neighborhood of $\cY$.
Let $\cY^{(2)}_{\lambda}$ be the corresponding space of $\lambda$-connections for $\lambda^2 = 0$ (explained in \cite{A04} in the case of connections without singularities).
Then, we will claim that for any function $f \in H^0(\cY^{(2)}_{\lambda}, \cO)$, its restriction to $\cY$ is constant.

This can be proven in the following way.
Denote by $S_m \subset S$ the open subset of the Hitchin base corresponding to smooth spectral curves.
Extended results of \cite{A04} provide an explicit description of an open subset of $\cY^{(2)}_{\lambda}$.
The description considers a particular moduli space of line bundles with connection $\cY^{\natural}_{sm}$, that is an affine bundle over $\cY_{sm} = \cY \times_{S} S_m$
and then defines a map
\begin{equation}
 \cY^{\natural}_{sm} \otimes \Spec(\bC[\lambda]/(\lambda^2)) \ra \cY^{(2)}_{\lambda}.
\end{equation}
When $\lambda = 0$, the map is simply the projection
\begin{equation}
 \cY^{\natural}_{sm} \ra \cY_{sm} \subset \cY.
\end{equation}
Moreover, there is an explicit description of the foliation that is tangent to the fibers of this map.

Now a regular function on $\cY^{(2)}_{\lambda}$ would give rise to a regular function on $\cY^{\natural}_{sm} \otimes \Spec(\bC[\lambda]/(\lambda^2))$ that is constant along this foliation.
However, any regular function on $\cY^{(2)}_{\lambda}$ comes from a regular function on $S_m \times \Spec(\bC[\lambda]/(\lambda^2))$.
The result follows from the fact that the foliation does not respect the projection $\cY^{\natural}_{sm} \ra S_m$, or the derivative of this foliation is transversal to the fiber of the projection.
\end{proof}

\begin{Rem}
\rm{The proof of Corollary \zcref{proof_of_main2} works for arbitrary $n \geq 5$ if we assume the statements corresponding to Proposition \zcref{cohomology_of_compactification} and Lemma \zcref{canonical_bundle_over_M_H}.}
\end{Rem}

\section{Orthogonality, conjectural compactified Radon transform, and Geometric Langlands Correspondence}\label{Section_GLC-like}
In this section, we will explain one way to extend Arinkin's results \cite{A01} via the Radon transform.
We learned the ideas explained here from D. Arinkin \cite{A}.
\subsection{Orthogonality}\label{section_orthogonality}
As explained in \cite[Section 3]{LS}, the moduli space $P_{-1}$ of indecomposable quasi-parabolic ${\mathfrak{s}l}_2$-bundles of degree $-1$ on $(\bP^1, D)$, where $D := t_1 + \cdots + t_n$, is a non-separated scheme which contains the projective space $\bP^{n-3}_{\bf b}$.
On the other hand, let us consider a new non-separated scheme that contains $(\bP^{n-3}_{\bf b})^{\vee} = \bP^{n-3}_{\bf a}$.
 Denote by $Z_i \subset \bP^{n-3}_{\bf a} = \Sym^{n-3}(\bP^1)$ the hyperplane of sections vanishing at $t_i \in \bP^1$ $(i = 1, \dots, n)$.
 
 \begin{Def}
  Let $P^{\vee}$ be the non-separated scheme obtained by gluing together two copies of $\bP^{n-3}_{\bf a}$ by the identity map over the open subset $U := \bP^{n-3}_{\bf a} \setminus \cup_{i = 1}^nZ_i$.
 \end{Def}
Note that $P_{-1} \simeq P^{\vee}$ in the $n=4$ case.
Let $Z_i^{\pm}$ be the pre-images of $Z_i \subset \bP^{n-3}_{\bf a}$ along $p \colon P^{\vee} \ra \bP^{n-3}_{\ba}$, and
$\bnu := \sum_{i = 1}^n \nu_i ([Z_i^+] - [Z_i^-]) \in \mathrm{div}(P^{\vee}) \otimes_{\bZ} \bC$, where $\mathrm{div}(P^{\vee})$ is the group of divisors on $P^{\vee}$.
Let $D_{\bnu}$ denote the ring of twisted differential operators (TDO) corresponding to $\bnu$ over $P^{\vee}$.

For any connection $\mathbb{L} = (L, \nabla, \varphi) \in \mathcal{M}$, its symmetric product $\Sym^{n-3}(\mathbb{L})$ gives a connection on $\mathbb{P}_{\bf a}^{n-3}$.
More precisely, it is the symmetric part of the push-forward of $\mathbb{L}^{\boxtimes n-3}$ along the map $\Sym \colon (\bP^1)^{n-3} \ra \bP^{n-3}_{\ba}$, that is, $\Sym^{n-3}(\mathbb{L}) := (\Sym_*(\mathbb{L}^{\boxtimes (n-3)}))^{\mathfrak{S}_{n-3}}$. 
This connection has singularities along the divisors $Z_i$ $(i = 1, \dots, n)$, as well as along the discriminant divisor $\Delta \subset \bP^{n-3}_{\bf a}$.
The divisors $Z_i$ cross normally, and the singularity along $Z_i$ has residue with eigenvalues $\{\pm \nu_i \}$, each with multiplicity $2^{n-4}$.

Let us construct the $D_{\bnu}$-module $ j_{!*}(\Sym^{n-3}(\mathbb{L})|_U)$ with $j \colon U := \bP^{n-3}_{\bf a} \setminus \cup_{i = 1}^n Z_i \hookrightarrow P^{\vee}$.
This construction still makes sense for a family of connections.
Let us apply it to the universal family of connections, and get a $\cM$-family $\xi_{\bnu}$ of $D_{\bnu}$-modules over $\cM \times P^{\vee}$. 
For ${\bm x} \in P^{\vee}$, denote by $(\xi_{\bnu})_{\bm x}$ the restriction of $\xi_{\bnu}$ to $\cM \times \{ {\bm x} \}$. 
Theorem \zcref{main} and its proof imply that $\xi_{\bnu}$ satisfies the orthogonality property over general points as follows:

\begin{Thm}\label{generic_orthogonal_property}
Suppose $n=5$, and ${\bm  x}, {\bm y} \in P^{\vee} \setminus (\cup_{i = 1}^5 Z_i^{\pm} \cup \Delta)$. Then
\begin{equation*}
 H^i(\cM, (\xi_{\bnu})_{\bm x} \otimes (\xi_{\bnu})_{\bm y}) = 0
\end{equation*}
for any ${\bm x} \neq {\bm y}$, $i \geq 0$.
\end{Thm}
\begin{proof}
If ${\bm x}, {\bm y} \in \bP^2_{\bf a} \setminus (\cup_{i = 1}^5 Z_i \cup \Delta)$ can be written as ${\bm x} = [x_1, x_2], {\bm y} = [ x_3, x_4]$ with $x_i \in \bP^1$ and $x_i \neq x_j$ for $i \neq j$, then the statement follows from Theorem \zcref{main}.

Suppose ${\bm x} = [x_1, x_2], {\bm y} = [ x_1, x_3]$.
Then, we need to check
$H^i(\cM, \xi_{x_1}^{\otimes 2} \otimes \xi_{x_2} \otimes \xi_{x_3}) = H^i(\cM, \xi_{x_2} \otimes \xi_{x_3}) \oplus H^i(\cM, \Sym^2(\xi_{x_1}) \otimes \xi_{x_2} \otimes \xi_{x_3}) = 0$.
But we can show $H^i(\cM, \xi_{x_2} \otimes \xi_{x_3}) = 0$ and $H^i(\cM, \Sym^2(\xi_{x_1}) \otimes \xi_{x_2} \otimes \xi_{x_3}) = 0$ in the same way as the proof of Theorem \zcref{main}.
\end{proof}

Moreover, it is predicted that for $n \geq 4$, the similar statement of Theorem \zcref{main} is also true with $2(n-3)$ points on $\bP^1$:
For $x \in \bP^1$ let $\xi_x$ be the bundle on $\overline{\cM}$ whose fiber at $(L, \nabla, \varphi; \e \in E)$ is $L_x$.
\begin{Conj}\label{generalization_of_main_thm}
Suppose $n\geq 4$, $x_1,\dots, x_{2(n-3)} \in \bP^1$ and $x_i \neq x_j$ for $i \neq j$. Then
\begin{equation*}
 H^i(\cM, \xi_{x_1}\otimes \cdots \otimes \xi_{x_{2(n-3)}}) = 0
\end{equation*}
 for any $i \geq 0$.
\end{Conj}
The $n=4$ case corresponds to \cite[Theorem 2 (i)]{A01}, and the $n=5$ case is our Theorem \zcref{main}.

In the subsequent sections, we will explain the meaning of Theorem \zcref{generic_orthogonal_property} from the viewpoint of geometric Langlands correspondence.

\subsection{Tamely ramified Geometric Langlands Correspondence}
We briefly recall the form of the tamely ramified geometric Langlands correspondence relevant to this paper.
For more detail, see \cite[Section 2]{A01}, \cite[Section 9]{AF12} and \cite[Appendix A]{DP}.

Let $C$ be a Riemann surface and $\Par_C = t_1 + \cdots + t_n$ be a reduced divisor on $C$.
In the ramified setting, the spectral side consists of logarithmic ${}^LG$-connections with prescribed local exponents along $\Par_C$, while the automorphic side consists of twisted $D$-modules on the moduli stack of quasi-parabolic $G$-bundles on $(C, \Par_C)$. 

We specialize to $G = PGL_2$, ${}^LG = SL_2$, and $(C, \Par_C) = (\bP^1, D := t_1 + \cdots + t_n)$.
In this case, the spectral side is the moduli stack $\cM$ of $\bnu$-${\mathfrak sl}_2$-parabolic connections on $(\bP^1, D)$.
On the automorphic side, the relevant component is ${\mathcal Bun}_{qp}^{odd}(PGL_2)$, the odd component of the moduli stack of quasi-parabolic $PGL_2$-bundles.
Its indecomposable part ${\mathcal Bun}_{qp}^{0, odd}(PGL_2)$ is identified with $P_{-1}$, the coarse moduli space of indecomposable quasi-parabolic ${\mathfrak sl}_2$-bundles of degree $-1$.

Denote by $\xi_i$ the invertible sheaf on ${\mathcal Bun}_{qp}^{0, odd}(PGL_2)$ whose fiber over $(\cF, \bl = \{ l_1, \dots, l_n \})$ is $l_i$ $(i = 1, \dots, n)$.
The local exponents $\bnu$ determine the element $\Sigma_i \nu_i [\xi_i] \in \Pic ({\mathcal Bun}_{qp}^{0, odd}(PGL_2)) \otimes_{\bZ} \bC$, where $[\xi_i] \in \Pic({\mathcal Bun}_{qp}^{0, odd}(PGL_2))$ is the isomorphic class of $\xi_i$.
It is predicted that there exists the so-called {\it{Okamoto map}} from $\Pic({\mathcal Bun}_{qp}^{0, odd}(PGL_2)) \otimes_{\mathbb{Z}}\bC$ to $\Pic(P_{-1})\otimes_{\mathbb{Z}}\bC$ (see  \cite[Section 7.6, Appendix A]{DP}).
Denote by $\Oka(\bnu)$ the image of $\Sigma_i \nu_i [\xi_i]$, and let $D(P_{-1})_{\bnu}$ be the TDO ring over $P_{-1}$ corresponding to $\Oka(\bnu)$.
The following conjecture is a version of {\it{tamely ramified Geometric Langlands Correspondence}} interpreted by D. Arinkin.

\begin{Conj}[{\cite[Section 2]{A01}}]\label{GLC}
Connected component of the derived category of quasi-coherent sheaves on $\cM$ is equivalent to the derived category of $D(P_{-1})_{\bnu}$-modules on $P_{-1}$:
 \begin{equation*}
  L \colon  \mathcal{D}_{qc}(\mathcal{M})^- \xrightarrow{\sim} \mathcal{D}(P_{-1}, D(P_{-1})_{\bnu}).
 \end{equation*}
 Here, $\cF \in \cD_{qc}(\cM)^{\pm}$ if and only if $-1 \in \mu_2$ acts on $H^i(\cF)$ as $\pm1$ for any $i$.
\end{Conj}
Arinkin proved the $n=4$ case in \cite{A01}.
In the next subsection, we replace the right-hand side by a category of twisted $D$-modules on the dual non-separated scheme $P^{\vee}$ via a conjectural compactified Radon transform.

\subsection{Radon Transform}
Let us recall the homogeneous Fourier transform, which is also known as the Radon transform.

Let $S$ be a projective space of dimension $d$, and let $D_{S, \chi}$ be a TDO ring on $S$ corresponding to a non-integral twist.
That is, if $S = \bP(V)$ for a $(d + 1)$-dimensional space $V$, then the category of $D_{S, \chi}$-modules is identified with the category of $D_V$-modules which is transformed by a fixed non-trivial character sheaf $\chi$ under the action of $\mathbb{G}_m$ by dilations.

Denote by $S^{\vee}$ the dual projective space, and by $D_{S^{\vee}, \chi^{-1}}$ the TDO ring on $S^{\vee}$ corresponding to the opposite twist.
Thus, the category of $D_{S^{\vee}, \chi^{-1}}$-modules is identified with the category of $D_{V^{\vee}}$-modules which is translated by $\chi^{-1}$ under the dilation action of $\mathbb{G}_m$.

The Radon transform is the exact equivalence
\begin{equation}\label{radon}
 R \colon \cD(S, D_{S, \chi}) \xrightarrow{\sim} \cD(S^{\vee}, D_{S^{\vee}, \chi^{-1}}),
\end{equation}
which descends from the Fourier transform between the corresponding categories of twisted equivariant $D$-modules on the corresponding vector spaces.

The product $S \times S^{\vee}$ carries a Poisson form which degenerates along the incidence variety $\Sigma \subset S \times S^{\vee}$.
Set $T := S \times S^{\vee} \setminus \Sigma$.
Then, the two projection $T \ra S$ and $T \ra S^{\vee}$ identify $T$ with the twisted cotangent bundles on $S$ and on $S^{\vee}$.
So, $T$ is a symplectic variety.
In our case, $T$ corresponds to $M^{{\bm w}_0}_{-1}$ with $S = \bP^{n-3}_{\bf b}$, $S^{\vee} = \bP^{n-3}_{\bf a}$ and $\Sigma \simeq M^{{\bm w}_0}_H$ (see \zcref{LSreview}).

Let us extend this picture to the whole $\widehat{M}'$, the coarse moduli space of $\bnu'$-$\mathfrak{sl}_2$-parabolic connections of degree $-1$ with a cyclic vector (see Definition \zcref{with_cyc_vec_odd}).
In this case, we replace $S$ (respectively, $S^{\vee}$) in \zcref{radon} with $P_{-1}$ (respectively, $P^{\vee}$).
This modification leads to the following conjecture:

\begin{Conj}[Partially Compactified Radon Transform]\label{Compactified_Radon}
The derived category of $D(P_{-1})_{\bnu}$-modules on $P_{-1}$ is equivalent to the derived category of $D_{\bnu}$-modules on $P^{\vee}$:
 \begin{equation*}
R' \colon   \mathcal{D}(P_{-1}, D(P_{-1})_{\bnu}) \xrightarrow{\sim} \mathcal{D}(P^{\vee}, D_{\bnu}).
 \end{equation*}
\end{Conj}
In the $n=4$ case, this equivalence is mentioned in \cite[Remark in Section 2]{A01}. 
Conjecture \zcref{Compactified_Radon} implies that Conjecture \zcref{GLC} is equivalent to:

\begin{Conj}\label{GLC-like}
Connected component of the derived category of quasi-coherent sheaves on $\cM$ is equivalent to the derived category of $D_{\bnu}$-modules on $P^{\vee}$:
 \begin{equation*}
  L' \colon  \mathcal{D}_{qc}(\mathcal{M})^- \xrightarrow{\sim} \mathcal{D}(P^{\vee}, D_{\bnu}).
 \end{equation*}
\end{Conj}

To establish the geometric Langlands correspondence, for each connection $\mathbb{L} \in \mathcal{M}$, we need to construct a $D(P_{-1})_{\bnu}$-module $\Aut_{\mathbb{L}}$ that satisfies the Hecke eigensheaf property.
It is predicted that $\Aut_{\mathbb{L}}$ is irreducible. Therefore, it suffices to describe its restriction to any open set.
$\Aut_{\mathbb{L}}$ is then recovered as the minimal extension from this open set.

Let $\iota : \bP^{n-3}_{\bf b} \hookrightarrow P_{-1}$ be a natural map.
Then, it is expected by D. Arinkin \cite{A} that $\iota^*(\Aut_{\mathbb{L}})$ is the $D(P_{-1})_{\bnu}|_{\bP^{n-3}_{\bf b}}$-module obtained by the Radon transform from $\Sym^{n-3}(\mathbb{L})$.

Therefore, let us consider $L'$ as a Fourier-Mukai transform with the kernel $\xi_{\bnu}$ constructed in \zcref{section_orthogonality}.
Corollary \zcref{proof_of_main2} and Theorem \zcref{generic_orthogonal_property} support the orthogonality property of $\xi_{\bnu}$ as an orthogonal $P^{\vee}$-family of $\cO_{\cM}$-modules in the $n = 5$ case.
In future work, we will check the orthogonality property along the remaining locus and prove the categorical equivalence (Conjecture \zcref{GLC-like}) in this case.

\section*{AI Declaration}
The manuscript was completed in October 2025 without the use of AI tools.

\bibliographystyle{alpha}
\bibliography{references_cohomology.bib}

\Address

\end{document}